\documentclass[11pt]{amsart}
\usepackage{amssymb, latexsym}
\theoremstyle{plain}
\newtheorem{theorem}{Theorem}
\newtheorem{corollary}{Corollary}

\newtheorem{proposition}{Proposition}

\newtheorem*{1'}{Theorem 1-Bessel}
\newtheorem*{P2'}{Proposition 2-Bessel}
\newtheorem*{P3'}{Proposition 3-Bessel}
\newtheorem*{P4'}{Proposition 4-Bessel}
\newtheorem*{C1'}{Corollary 1-Bessel}

\newtheorem*{2'}{Theorem 2-Bessel}
\newtheorem*{3'}{Theorem 3-Bessel}

\theoremstyle{remark}

\newtheorem*{Remark 1}{Remark 1}
\newtheorem*{Remark 2}{Remark 2}
\newtheorem*{Remark 3}{Remark 3}
\newtheorem*{Remark 4}{Remark 4}

\numberwithin{equation}{section}
\renewcommand{\baselinestretch}{1.22}

\begin{document}

\title [Probability of failure in diffusive search with resetting]
{Large time probability of failure in diffusive search with resetting for a random target in $\mathbb{R}^d$--a functional analytic approach}

\author{Ross G. Pinsky}

%\noindent  pinsky@math.technion.ac.il\ \ \ \ tel: 972-4-829-4083\ \ \  fax: 972-4-829-3388

\address{Department of Mathematics\\
Technion---Israel Institute of Technology\\
Haifa, 32000\\ Israel}
\email{ pinsky@math.technion.ac.il}

\urladdr{http://www.math.technion.ac.il/~pinsky/}

\subjclass[2010]{60J60} \keywords{random target, diffusive search, resetting,  principal eigenvalue, large time asymptotic behavior}
\date{}

\begin{abstract}
We consider a stochastic search model with resetting for an unknown stationary target $a\in\mathbb{R}^d,\ d\ge1$, with known distribution $\mu$. The searcher begins at the origin and performs Brownian motion
with diffusion coefficient $D$. The searcher is also equipped with an exponential clock with rate $r>0$, so that if
 it has  failed to locate the target by the time the clock rings, then its position is reset to the origin and it   continues its search anew from there. In dimension one, the target is considered located when the process
hits the point $a$, while in dimensions two and higher, one chooses an $\epsilon_0>0$ and the target is considered located
when the process hits the $\epsilon_0$-ball centered at $a$.
  Denote the position of the searcher at time $t$  by $X(t)$,
  let $\tau_a$ denote the time that a target at $a$ is located, and let $P^{d;(r,0)}_0$
  denote probabilities for
the process starting from 0. Taking a
functional analytic  point of view, and
  using the generator of the Markovian search process and its adjoint,
  we obtain
  precise estimates, with control on the dependence on   $a$, for the asymptotic behavior of $P^{d;(r,0)}_0(\tau_a>t)$ for large time,
  and then use this to obtain large time estimates on
  $\int_{\mathbb{R}^d}P^{d;(r,0)}_0(\tau_a>t)d\mu(a)$,
  the probability that
the searcher has failed up to time $t$ to locate the random target, for a variety of families of target distributions $\mu$.
Specifically,
for $B,l>0$ and $d\in\mathbb{N}$, let $\mu^{(d)}_{B,l}\in\mathcal{P}(\mathbb{R}^d)$ denote
any  target distribution   with density $\mu_{B,l}^{(d)}(a)$ that  satisfies
$$
\lim_{|a|\to\infty}\frac{\log\mu_{B,l}^{(d)}(a)}{|a|^l}=-B.
$$
Then we prove that
$$
\lim_{t\to\infty}\frac1{(\log t)^l}\log\int_{\mathbb{R}^d}  P_0^{d;(r,0)}(\tau_a>t)\mu_{B,l}^{(d)}(da)=-B(\frac D{2r})^\frac l2.
$$
The result is independent of the dimension.
In particular, for example, if the target distribution is
a centered Gaussian  of any dimension with variance $\sigma^2$, then
for any $\delta>0$, the probability of not locating the target by time $t$ falls in the interval
$\big(e^{-(1+\delta)\frac{D}{4r\sigma^2}(\log t)^2}, e^{-(1-\delta)\frac{D}{4r\sigma^2}(\log t)^2}\big)$, for sufficiently large $t$.

\end{abstract}

\maketitle
\section{Introduction and Statement of Results}\label{intro}
\renewcommand{\baselinestretch}{1.33}
The use of resetting in  search problems is
 a  common phenomenon in various contexts.  For example, in everyday life,  one might be searching for some target, such as a face in a crowd or a misplaced object. After having searched unsuccessfully for a while,
  there is   a tendency to return to the starting point  and begin the search anew.
Other contexts where search problems  frequently involve resetting include  animal foraging \cite{BC,VdRS},  proteins searching for target sites  on DNA molecules \cite{BWv,CBM,GMKC} and internet search algorithms.

Over the past decade or so, a variety of  stochastic processes with resetting have attracted much attention, mainly in the physics literature. See \cite{EMS} for a rather comprehensive, recent overview.  Prominent among such processes is the  diffusive search process with resetting, the process we consider in this paper.
Consider a  random  stationary target  $a\in \mathbb{R}^d$
with  known distribution $\mu$, and consider  a searcher who sets off from the origin, and performs  $d$-dimensional Brownian motion with diffusion coefficient $D$.
The searcher is also equipped with   an exponential clock with  rate $r$, so that
if it has  failed to locate the target by the time the clock rings, then its position is reset to the origin and it continues  its search anew from there.
In dimension one, the target is considered ``located'' when the process
hits the point $a$, while in dimensions two and higher, one chooses an $\epsilon_0>0$ and the target is considered ``located''
when the process hits the $\epsilon_0$-ball centered at $a$.
One may be interested in several statistics, the most important ones being  the expected time to locate the target and the probability of failing to reach the target after a large time.
  See, for example, \cite{EM1, EM2,EM3, EMM, FBPCM, P, NG, dMMT, RR} for a sampling of articles on this model and related ones.

The objective of this paper is to give a rigorous analysis of
the latter of these two statistics, from a functional analytic  point of view,
  using the generator of the Markovian search process and its adjoint.
However,  we begin with some comments concerning the first of these statistics.
 Without the resetting, the expected time to locate the target at  any fixed $a\in\mathbb{R}^d-\{0\}$  is infinite \cite{Pin95}.
With the resetting,
%Under very general conditions on the    rate of the resetting time, and in particular, when the rate is constant,
the expected time to locate the target at $a\in \mathbb{R}^d$ is finite. In dimension one it is  given
by $\frac{e^{\sqrt{\frac{2r}D}\thinspace |a|}-1}r,\ a\in\mathbb{R}$, \cite{EM2} while in dimensions $d\ge2$ it is given explicitly in terms of the modified Bessel function of the second kind, $K_{\frac{d-2}2}$ \cite{EM3}.
From the above formula in one-dimension, the expected time to locate the random target
is $\int_{-\infty}^\infty  \frac{e^{\sqrt{\frac{2r}D}\thinspace |a|}-1}r\mu(da)$.
In particular, in order for this expected time  to be finite, the target distribution $\mu$ must possess  some exponential moments.
A similar phenomenon holds in higher dimensions.
In \cite{Pin20}, a spatially dependent exponential resetting rate was considered in the one-dimensional case, and it was shown that for any distribution $\mu$ with finite $l$th moment, for some $l>2$, one can choose
a spatially dependent resetting rate so that the expected time to locate the random target is finite.

In this paper we consider a constant resetting rate $r$.
Before discussing our results concerning  the large time probability that the searcher fails to locate the target, we
 give a more formal  mathematical definition of the model.
The process $X(t)$ on $\mathbb{R}^d$ is defined as follows.
The process starts from $0\in\mathbb{R}^d$ and  performs $d$-dimensional Brownian motion with diffusion coefficient $D$, until a random  clock rings.
This random clock has an exponential distribution with parameter $r$, so the probability that it has not rung by time $t$ is $e^{-rt}$.
When the clock rings, the process is instantaneously reset to its initial position 0, and continues its search  afresh with an independent resetting clock, and the above scenario is repeated, etc.
We define the process so that it is right-continuous.
Denote probabilities and expectations for the process starting at $x\in\mathbb{R}^d$ by $P^{d;(r,0)}_x$ and $E^{d;(r,0)}_x$ respectively.
The pair $(r,0)$ in the notation refers to the resetting rate $r$ and the resetting position $0$. (For the analysis in the multidimensional case, we will need to consider resetting to a point different than 0.)
From the above description, it follows that  $X(t)$ is a Markov process whose generator
$L^{d;(r,0)}$, restricted to appropriate functions $u$,    satisfies
\begin{equation}\label{generator}
L^{d;(r,0)}_du(x)=\frac D2\Delta u(x)+r\big(u(0)-u(x)\big).
\end{equation}
(See the proof of Proposition \ref{Ttcompact} and Proposition 2-Bessel for more details.)

Fix  $\epsilon_0>0$ once and for all.
Let
\begin{equation}\label{stoppingtime}
\tau_a=\begin{cases}\inf\{t\ge0:X(t)=a\},\ d=1;\\ \inf\{t\ge0:|X(t)-a|\le\epsilon_0\},\ d\ge 2\end{cases}
\end{equation}
denote the time at which a target at  $a\in\mathbb{R}^d$ is located.
In this paper, we
%use  rigorous   functional analytic and probabilistic methods to
study the  asymptotic behavior as $t\to\infty$ of  $P_0^{d;(r,0)}(\tau_a>t)$, the probability that the resetting process has not located a target at  $a$ by time $t$, and then   use this to  analyze
the asymptotic behavior as  $t\to\infty$ of  $\int_{-\infty}^\infty  P_0^{d;(r,0)}(\tau_a>t)\mu(da)$,  the probability that
the searcher has failed up to time $t$ to locate the random target, distributed according to $\mu\in\mathcal{P}(\mathbb{R}^d)$.
The asymptotic behavior  of  $P_0^{d;(r,0)}(\tau_a>t)$ has already been investigated in \cite{EM1} for the one-dimensional case and in \cite{EM3} for the multi-dimensional case,  using the method of inverse Laplace transforms. The mathematics there is a bit informal.
Using our functional analytic approach,
the basic asymptotic behavior we obtain is the same   as in those papers, however
the form in which we obtain it gives us explicit control over the dependence of this behavior on  $a$, in contrast to the state of affairs in the above-mentioned papers, as far as this author can tell.
We elaborate on this more in the next paragraph.
 This control is crucial for the next step, which is the main point of the paper,  namely the analysis of
$\int_{-\infty}^\infty  P_0^{d;(r,0)}(\tau_a>t)\mu(da)$.
In addition, the form in which we obtain our estimate
on $P_0^{d;(r,0)}(\tau_a>t)$
  allows for greater understanding of the underlying probabilistic mechanisms at work.
 Furthermore, we identify explicitly a number of spectral theoretic quantities, such as the principal eigenfunctions
 of the operator and its adjoint, and this might be of some independent interest.
 The papers in the physics literature have not studied the asymptotic behavior of
  $\int_{-\infty}^\infty  P_0^{d;(r,0)}(\tau_a>t)\mu(da)$,  the probability that
the searcher has failed to locate the random target by time $t$; thus, our work on this is entirely new.

An asymptotic formula of the form $P_0^{d;(r,0)}(\tau_a>t)\sim c(a,t)e^{-\lambda_0(r,0;a)t}$ is obtained
both in \cite{EM1,EM3} and in this paper, where $\lambda_0(r,0;a)$ satisfies a certain implicit
equation, which allows for its asymptotic analysis as $a\to\infty$.
In \cite{EM1,EM3}, $\lambda_0(r,0;a)$
arises from the inverse Laplace transform method, while in this paper, it arises as a certain principal eigenvalue.
However, the term $c(a,t)$ is not analyzed sufficiently for our needs  in \cite{EM1,EM3}.
In our paper, we obtain the term $c(a,t)$ explicitly in terms of an expectation involving the search process, and this allows us sufficient control over $c(a,t)$ in order to study
the asymptotic behavior of  $\int_{-\infty}^\infty  P_0^{d;(r,0)}(\tau_a>t)\mu(da)$ for certain families  of target distributions $\mu$.

Before stating the main results,
we describe a  side result which will follow readily from the results concerning  $P_0^{d;(r,0)}(\tau_a>t)$.
The  Brownian motion without resetting  corresponds to setting $r=0$; let  $P_x^{d;(0)}$ and $E_x^{d;(0)}$ denote
probabilities and expectations for the  Brownian motion without resetting starting from $x\in \mathbb{R}^d$.
As already noted, for fixed $a\in\mathbb{R}^d-\{0\}$, the expected time to locate a target at $a$ by a Brownian motion without resetting is infinite,
but for the Brownian motion with resetting it is finite.
 However,  the one-dimensional (two-dimensional) Brownian motion without resetting reaches distant points (the $\epsilon_0$-neighborhood of distant points) much more quickly than does one-dimensional (two-dimensional) Brownian motion with resetting.
 (Of  course, in three dimensions and higher, Brownian motion without resetting has a positive probability of never reaching the $\epsilon_0$-neighborhood of a point.)
 In the one-dimensional case without resetting,
using Brownian scaling (or alternatively, the reflection principle), one can readily show  that
$$
\lim_{t\to\infty}P_0^{1;(0)}(\tau_{a_t}>t)=
\begin{cases} 0,\ \text{if}\ \lim_{t\to\infty}\frac{|a_t|}{\sqrt t}=0;\\ 1,\ \text{if}\   \lim_{t\to\infty}\frac{|a_t|}{\sqrt t}=\infty.\end{cases}
$$
In the two-dimensional case without resetting,
we have the following result.
\begin{proposition}\label{2dnoreset}
\begin{equation}\label{2dnoresetresult}
\begin{aligned}
&\lim_{t\to\infty}P_0^{2;(0)}(\tau_{a_t}>t)=0,\ \text{if}\ \lim_{t\to\infty}\frac{|a_t|}{t^\delta}=0,\ \text{for all}\ \delta>0;\\
&\liminf_{t\to\infty}P_0^{2;(0)}(\tau_{a_t}>t)>0,\ \text{if}\ \liminf_{t\to\infty}\frac{|a_t|}{t^\delta}>0,\ \text{for some}\ \delta\in(0,\frac12];\\
&\lim_{t\to\infty}P_0^{2;(0)}(\tau_{a_t}>t)=1,\ \text{if}\ \lim_{t\to\infty}\frac{|a_t|}{t^\frac12}=\infty.
\end{aligned}
\end{equation}
\end{proposition}

On the other hand, we will prove the following result for the Brownian motion with resetting.
\begin{proposition}\label{speed}
For $d=1$,
\begin{equation}\label{speedequ}
\lim_{t\to\infty}P_0^{1;(r,0)}(\tau_{a_t}>t)=\begin{cases}0,\ \text{\rm if}\ \lim_{t\to\infty}(|a_t|-\sqrt{\frac D{2r}}\log t)=-\infty;\\ 1,\ \text{\rm if}\  \lim_{t\to\infty}(|a_t|-\sqrt{\frac D{2r}}\log t)=\infty.\end{cases}
\end{equation}
For $d\ge2$,
\begin{equation}\label{speedequd}
\begin{aligned}
&\lim_{t\to\infty}P_0^{d;(r,0)}(\tau_{a_t}>t)=0,\ \text{\rm if}\ \lim_{t\to\infty}(|a_t|-\sqrt{\frac D{2r}}\log t+\gamma\log\log t)=-\infty,\\
& \text{\rm for some}\ \gamma>\frac{d-1}2\sqrt{\frac D{2r}};\\
&\lim_{t\to\infty}P_0^{d;(r,0)}(\tau_{a_t}>t)=1,\ \text{\rm if}\ \lim_{t\to\infty}(|a_t|-\sqrt{\frac D{2r}}\log t+\frac{d-1}2\sqrt{\frac D{2r}}\log\log t)=\infty.
\end{aligned}
\end{equation}
\end{proposition}

\medskip

We now turn to the main results. We will be interested in the behavior of
the process with $D$ and $r$ fixed. In our notation, we suppress all dependence on $D$
(except in Corollary \ref{cor1} and  Corollary 1-Bessel, where the dependence of certain constants on $D$ is indicated),
but indicate the dependence
on $r$.
We begin by stating our central result, which  concerns
$\int_{-\infty}^\infty  P_0^{d;(r,0)}(\tau_a>t)\mu(da)$,  the probability that
the searcher has failed to locate the random target by time $t$.

For $B,l>0$ and $d\in\mathbb{N}$, let $\mu^{(d)}_{B,l}\in\mathcal{P}(\mathbb{R}^d)$ denote
any  target distribution   with density $\mu_{B,l}^{(d)}(a)$ that  satisfies
\begin{equation}\label{targetdist}
\lim_{|a|\to\infty}\frac{\log\mu_{B,l}^{(d)}(a)}{|a|^l}=-B.
\end{equation}

\begin{theorem}\label{randomtarget}
Let $\mu_{B,l}^{(d)}\in\mathcal{P}(\mathbb{R}^d)$ be a distribution  with density  satisfying \eqref{targetdist}.
Then $\int_{\mathbb{R}^d}  P_0^{d;(r,0)}(\tau_a>t)\mu_{B,l}^{(d)}(da)$, the probability that the searcher with resetting fails to locate the random target with distribution
$\mu_{B,l}^{(d)}$ by time $t$, satisfies
\begin{equation}\label{asympprobtarget}
\lim_{t\to\infty}\frac1{(\log t)^l}\log\int_{\mathbb{R}^d}  P_0^{d;(r,0)}(\tau_a>t)\mu_{B,l}^{(d)}(da)=-B(\frac D{2r})^\frac l2.
\end{equation}
\end{theorem}
\noindent \bf Remark.\rm\ Unlike all of the other results in this paper, the result in Theorem \ref{randomtarget} is
independent of the  dimension.

\medskip

\noindent \bf Example 1.\rm\ Consider  a target distribution of the form \eqref{targetdist} with $l=1$.
In particular, if $d=1$, this situation includes
 the   two-sided, symmetric exponential distributions, whose densities
are of the form $Be^{-B|x|}$, $B>0$.
 One has that for any $\delta>0$,  the probability of not locating the target by time $t$ falls in the
interval  $(t^{-(1+\delta)B\sqrt{\frac D{2r}}},t^{-(1-\delta)B\sqrt{\frac D{2r}}})$, for sufficiently large $t$.
\medskip

\bf\noindent Example 2.\rm\ Consider a centered Gaussian target distribution in any dimension, with variance $\sigma^2$. This distribution is of the form \eqref{targetdist} with $l=2$ and $B=\frac1{2\sigma^2}$.
For such a target distribution, for any $\delta>0$, the probability of not locating the target by time $t$ falls in the interval
$\big(e^{-(1+\delta)\frac{D}{4r\sigma^2}(\log t)^2}, e^{-(1-\delta)\frac{D}{4r\sigma^2}(\log t)^2}\big)$, for sufficiently large $t$.
\medskip

For the rest of the results, we need to treat separately  the one-dimensional and the multi-dimensional cases.
We begin with the one-dimensional case. We present a series of results which culminates in a formula for $P_0^{1;(r,0)}(\tau_a>t)$ of the form
$c(a,t)e^{-\lambda_0(r,0;a)t}$, where $\lambda_0(r,0;a)$ is a certain principal eigenvalue and $c(a,t)$ is given in terms of a certain conditional expectation (Theorem \ref{Main}), and a result which
estimates $c(a,t)$ for large $a$, uniformly in $t$ (Proposition \ref{alimitunift}).

For $a\neq0$,
let $T^{(r,0;a)}_t$ denote the semigroup corresponding to the Markov process $X(t)$ that is killed upon reaching $a$.
If $a>0$, then
\begin{equation}\label{Tt}
T^{(r,0;a)}_tf(x)=E^{1;(r,0)}_x(f(X(t)); \tau_a>t),\ x\in(-\infty,a], \ t\ge0,
\end{equation}
for bounded functions $f$ defined on $(-\infty,a)$.
For $a<0$, we have the corresponding formula with
$x\in[a,\infty)$.
From now on we will assume that $a>0$; of course all the results also hold for $a<0$, mutatis mutandis.
Let $[-\infty,a]$ denote the one-point compactification of $(-\infty,a]$, obtained by adding the point at $-\infty$,
and let $C_{0_a}\big([-\infty,a]\big)$
denote the space of continuous functions on $[-\infty,a]$ which vanish at $a$.
(Note that this space is equivalent to the space of continuous functions $u$ on $(-\infty,a)$ which satisfy
$\lim_{x\to-\infty}u(x)$ exists and $\lim_{x\to a}u(x)=0$.)
We will prove the following proposition. As usual, $C^2_b((-\infty,a))$ denotes the space
of functions defined on $(-\infty, a)$ which have two continuous and bounded derivatives.
\begin{proposition}\label{Ttcompact}
For $a,r>0$ and all $t>0$,  the semigroup operator $T^{(r,0;a)}_t$ is compact from
$C_{0_a}\big([-\infty,a]\big)$ to $C_{0_a}\big([-\infty,a]\big)$. Furthermore, its generator, which we denote
by $L^{(r,0;a)}$, is
an extension of the operator $L^{1;(r,0)}$ in \eqref{generator} defined
on $C^2_b((-\infty,a))\cap\{f:f,L^{1;(r,0)}f\in C_{0_a}\big([-\infty,a]\big)\}$.
\end{proposition}
From Proposition \ref{Ttcompact} it follows that the generator  $L^{(r,0;a)}$ has a compact resolvent
and consequently a principal eigenvalue, which we denote by $\lambda_0(r,0;a)$.
The following theorem and corollary concern this  principal eigenvalue and the corresponding principal eigenfunction.
\begin{theorem}\label{prinev}
Let $a>0$. The principal eigenvalue $\lambda_0(r,0;a)$ of the generator $L^{(r,0;a)}$
of the semigroup  $T^{(r,0;a)}_t$ is the unique solution $\lambda\in(0,r)$ of the equation
\begin{equation}\label{evequation}
\lambda=r\exp\big(-a\sqrt{\frac2D(r-\lambda)}\thinspace\big).
\end{equation}
A corresponding principal eigenfunction $u_{r,0;a}$ is given by
\begin{equation}\label{preigenfu}
u_{r,0;a}(x)=\frac{r}{r-\lambda_0(r,0;a)}(1-\exp\big(-\sqrt{\frac{2\big(r-\lambda_0(r,0;a)\big)}D}\thinspace(a-x)\big),\ x\le a.
\end{equation}
\end{theorem}
\medskip

\begin{corollary}\label{cor1}
\begin{equation}\label{2sidedevbdd}
re^{-\sqrt{\frac{2r}D}a}\le\lambda_0(r,0;a)\le re^{-c\sqrt{\frac{2r}D}a},
\end{equation}
where $c=c(r,a,D)\in(0,1)$ and $\lim_{a\to\infty} c(r,a,D)=1$.
\end{corollary}

For the statement and proof of Theorem \ref{Main} below we need to introduce the adjoint semigroup $\tilde{T_t}^{(r,0;a)}$ to the semigroup $T_t^{(r,0;a)}$. Since $T_t^{(r,0;a)}$ is  defined on the Banach space  $C_{0_a}([-\infty,a])$,
the adjoint $\tilde{T_t}^{(r,0;a)}$  operates on the dual space of bounded linear functions on $C_{0_a}([-\infty,a])$.
Since
$[-\infty,a]$ with the one-point compactification topology is compact,
this dual space is the space of finite signed measures on $[-\infty,a)$ \cite[p.28]{SS}.
Recall that a finite signed measure $\nu$ is of the form $\nu=\nu^+-\nu^-$, where $\nu^+,\nu^-$ are finite measures.
 (The reason the measures are on $[-\infty,a)$ instead of on $[-\infty, a]$ is that $f(a)=0$, for $f\in C_{0_a}([-\infty,a])$.)
Let $\nu$ be such a finite signed measure. We can write $\nu=c_+\nu^+-c_-\nu^-$,
where $\nu^+$ and $\nu^-$ are probability measures on $[-\infty,a)$ and $c_+,c_-\ge0$.
From \eqref{Tt}, it follows that
\begin{equation}\label{Ttstar}
\begin{aligned}
&\tilde T_t^{(r,0;a)}\nu(dy)=c_+\int_{-\infty}^a\nu^+(dx)P_x^{(r,0;a)}(X(t)\in dy;\tau_a>t)-\\
&c_-\int_{-\infty}^a\nu^-(dx)P_x^{(r,0;a)}(X(t)\in dy;\tau_a>t).
\end{aligned}
\end{equation}

Denote the generator of the adjoint semigroup by $\tilde{L}^{(r,0;a)}$. Of course, this operator has the same principal eigenvalue as does $L^{(r,0;a)}$.
 \begin{proposition}\label{adjoint}
The generator   $\tilde{L}^{(r,0;a)}$ of $\tilde{T_t}^{(r,0;a)}$
satisfies
\begin{equation}\label{adjointoper}
 \tilde{L}^{(r,0;a)}v(y)=\frac D2v''(y)-rv(y)+r\big(\int_{-\infty}^av(x)dx\big)\delta_0(y),\
\end{equation}
for $v$ satisfying $v\in C_{0_a}([-\infty,a])\cap C_b^2((-\infty, a))$ and $\int_{-\infty}^a|v(y)|dy<\infty$.
Furthermore, a principal eigenfunction $v_{r,0;a}$  corresponding to the principal eigenvalue $\lambda_0(r,0;a)$ is given by
\begin{equation}\label{adjointpreigenfu}
v_{r,0;a}(y)=\begin{cases} \exp\big(\sqrt{\frac{2(r-\lambda_0(r,0;a))}D}\thinspace y\big),\ y<0;\\
\frac{\exp\big(-\sqrt{\frac{2(r-\lambda_0(r,0;a))}D}\thinspace(y-a)\big)-
\exp\big(\sqrt{\frac{2(r-\lambda_0(r,0;a))}D}\thinspace(y-a)\big)}
{\exp\big(\sqrt{\frac{2(r-\lambda_0(r,0;a))}D}\thinspace a\big)-\exp\big(-\sqrt{\frac{2(r-\lambda_0(r,0;a))}D}\thinspace a\big)}\end{cases},\ 0\le y\le a.
\end{equation}
 \end{proposition}
 \medskip

\bf \noindent Remark.\rm\ The right hand side of \eqref{adjointoper} should be understood as
the signed measure whose absolutely continuous part has density
$\frac D2v''(y)-rv(y)$, and whose singular part is $r\big(\int_{-\infty}^av(y)dy\big)\delta_0$.

%The function $v$ is continuous but has a discontinuity in its derivative at $y=0$.

Here is our  result concerning the asymptotic behavior of $P_0^{1;(r,0)}(\tau_a>t)$.
\begin{theorem}\label{Main}
Let $a>0$. Then
\begin{equation}\label{asympprob}
P_0^{1;(r,0)}(\tau_a>t)=\frac1{E_0^{1;(r,0)}(u_{r,0;a}(X(t))|\tau_a>t)}e^{-\lambda_0(r,0;a)\thinspace t},
\end{equation}
where $u_{r,0;a}$ is as in \eqref{preigenfu},
Furthermore,
\begin{equation}\label{M}
\begin{aligned}
&\lim_{t\to\infty}E_0^{1;(r,0)}(u_{r,0;a}(X(t))|\tau_a>t)=
\frac{\int_{-\infty}^au_{r,0;a}(x)v_{r,0;a}(x)dx}{\int_{-\infty}^av_{r,0;a}(x)dx}=\\
&\thinspace\frac{2e^{qa}-2-qa}{2e^{qa}(1-\frac{\lambda(r,0;a)}r)^2}, \ \text{with}\  q=\sqrt{\frac{2(r-\lambda(r,0;a))}D},
\end{aligned}
\end{equation}
where    $v_{r,0;a}$ is as in \eqref{adjointpreigenfu}.
Thus, for fixed $a$,
\begin{equation}\label{tasymp}
\begin{aligned}
&P_0^{1;(r,0)}(\tau_a>t)\sim
\frac{2e^{qa}(1-\frac{\lambda(r,0;a)}r)^2}{2e^{qa}-2-qa}
e^{-\lambda_0(r,0;a)\thinspace t}, \text{as}\ t\to\infty,\\
&\text{with}\  q=\sqrt{\frac{2(r-\lambda(r,0;a))}D}.
\end{aligned}
\end{equation}
\end{theorem}

The following proposition concerns the coefficient multiplying the exponential term
in \eqref{asympprob}. It will be needed for the proof of Theorem \ref{randomtarget} as well as for the proof of Proposition \ref{speed}.
\begin{proposition}\label{alimitunift}
\begin{equation}\label{atoinf}
\lim_{a\to\infty}E_0^{1;(r,0)}(u_{r,0;a}(X(t))|\tau_a>t)=1,\ \text{uniformly over}\ t\in(0,\infty),
\end{equation}
where $u_{r,0;a}$ is as in \eqref{preigenfu}.
\end{proposition}

We now turn to the multi-dimensional case. Recall the definition of $\tau_a$ from \eqref{stoppingtime}.
We  make a construction to reduce the study of
$P_0^{d;(r,0)}(\tau_a>t)$
 to a one-dimensional problem. Instead  of having the target at $a\in\mathbb{R}^d$ and having the resetting bring the process to $0\in\mathbb{R}^d$, we  consider the target to be at 0 and have the resetting bring the process to $a$.
If we denote this new process by $\hat X(t)$ and denote probabilities by
$\hat P_x^{d;(r;a)}$, then clearly
\begin{equation}\label{sameprob1}
P_0^{d;(r,0)}(\tau_a>t)=\hat P_a^{d;(r,a)}(\hat\tau_0>t),
\end{equation}
where, consistent with the notation in \eqref{stoppingtime},
$$
\hat \tau_0=\inf\{t\ge0:|\hat X(t)|\le \epsilon_0\}.
$$
Now let $Y(t)=|\hat X(t)|$.
Then $Y(t)$ is the radial part of a $d$-dimensional Brownian motion with diffusion coefficient $D$,
and it is reset at rate $r$ to $|a|$. That is, $Y(t)$ is a Bessel process with resetting, of order $d$ with diffusion coefficient $D$.
% so its generator is
%$\frac{D}2\frac{d^2}{dx^2}+D\frac{d-1}{2x}\frac{d}{dx}$ on $(0,\infty)$.
%   the resetting of $\hat X(t)$ brings $Y(t)$ to
%$|a|$ and $\hat\tau_{0;a}$ is the first time $Y(t)\le \epsilon_0$. Let
Let
$$
\tau^{(\text{Y})}_{\epsilon_0}=\inf\{t\ge0:Y(t)=\epsilon_0\}.
$$
Denote probabilities and expectations for  $Y(t)$ starting at $x>\epsilon_0$ and with resetting to $A\in(\epsilon_0,\infty)$ at rate $r$  by $\mathcal{P}_x^{(r,A)}$ and $\mathcal{E}_x^{(r,A)}$.
Then clearly,
\begin{equation}\label{sameprob2}
\hat P_a^{d;(r,a)}(\hat\tau_0>t)=\mathcal{P}_{|a|}^{(r,|a|)}(\tau_{\epsilon_0}^{(Y)}>t),\ |a|>\epsilon_0.
\end{equation}
From \eqref{sameprob1} and \eqref{sameprob2}, it
follows that for the analysis of $P_0^{d;(r,0)}(\tau_a>t)$,  it suffices to study
$\mathcal{P}_{|a|}^{(r,|a|)}(\tau_{\epsilon_0}^{(Y)}>t)$.

We now present the analogs of Proposition \ref{Ttcompact}, Theorem \ref{prinev}, Corollary \ref{cor1},
Proposition \ref{adjoint}, Theorem \ref{Main} and Proposition \ref{alimitunift}
in the context of the above Bessel process with resetting.
We use the same labelling and numbering of theorems, propositions and the corollary as was used in the one-dimensional case, but suffix each of these with ``Bessel''.

The generator of the Bessel process of order $d$ with diffusion coefficient $D$ is
$\frac D2\frac{d^2}{dx^2}+D\frac{d-1}{2x}\frac d{dx}$. Define the operator $\mathcal{L}^{(r,A)}$ by
\begin{equation}\label{genmult}
\mathcal{L}^{(r,A)}u(x)=\frac D2u''(x)+D\frac{d-1}{2x}u'(x)+r\big(u(A)-u(x)\big).
\end{equation}
For $A>0$,
let $\mathcal{T}^{(r,A;\epsilon_0)}_t$ denote the semigroup corresponding to the Markov process $Y(\cdot)$ with resetting to $A$ at rate $r$, and which is killed upon  reaching $\epsilon_0$.
Then
\begin{equation}\label{Ttmult}
\mathcal{T}^{(r,A;\epsilon_0)}_tf(x)=\mathcal{E}^{(r,A)}_x(f(Y(t)); \tau^{(Y)}_{\epsilon_0}>t),\ x\in[\epsilon_0,\infty), \ t\ge0.
\end{equation}
Let $[\epsilon_0,\infty]$ denote the one-point compactification of $[\epsilon_0,\infty)$, obtained by adding the point at $\infty$,
and let $C_{0_{\epsilon_0}}\big([\epsilon_0,\infty]\big)$
denote the space of continuous functions on $[\epsilon_0,\infty]$ which vanish at $\epsilon_0$.
%\begin{proposition}\label{Ttcompactmult}
\begin{P2'}
For $A,r>0$ and all $t>0$,  the semigroup operator $\mathcal{T}^{(r,A;\epsilon_0)}_t$ is compact from
$C_{0_{\epsilon_0}}\big([\epsilon_0,\infty]\big)$ to $C_{0_{\epsilon_0}}\big([\epsilon_0,\infty]\big)$. Furthermore, its generator, which we denote
by $\mathcal{L}^{(r,A;\epsilon_0)}$, is
an extension of the operator $\mathcal{L}^{(r;A)}$ in \eqref{genmult} defined
on $C^2_b((\epsilon_0,\infty))\cap\{f:f,f',f''\in C_{0_{\epsilon_0}}\big([\epsilon_0,\infty]\big)\}$.
\end{P2'}

From Proposition 2-Bessel, it follows that the generator  $\mathcal{L}^{(r,A;\epsilon_0)}$ has a compact resolvent
and consequently a principal eigenvalue, which we denote by $\lambda_0(r,A;\epsilon_0)$.
The following theorem and corollary concern this  principal eigenvalue and the corresponding principal eigenfunction.
In the sequel, $K_\nu$ denotes the modified Bessel function of the second kind of order $\nu$.
This function  decays exponentially at $\infty$ \cite{AS,W}.
\begin{2'}
The principal eigenvalue $\lambda_0(r,A;\epsilon_0)$ of the generator $\mathcal{L}^{(r,A;\epsilon_0)}$
of the semigroup  $\mathcal{T}^{(r,A;\epsilon_0)}_t$ is the unique solution $\lambda\in(0,r)$ of the equation
\begin{equation}\label{evequationmult}
\lambda=r\big(\frac {A}{\epsilon_0}\big)^{\frac{2-d}2}\frac{K_{\frac{d-2}2}(\sqrt{(r-\lambda)\frac2D}\thinspace A)}{K_{\frac{d-2}2}(\sqrt{(r-\lambda)\frac2D}\thinspace\epsilon_0)}.
%\exp\big(-a\sqrt{\frac2D(r-\lambda)}\thinspace\big).
\end{equation}
A corresponding principal eigenfunction $\mathcal{U}_{r,A;\epsilon_0}$ is given by
\begin{equation}\label{preigenfumult}
\mathcal{U}_{r,A;\epsilon_0}(x)=\frac{r}{r-\lambda_0(r,A;\epsilon_0)}\Big(1-\big(\frac {x}{\epsilon_0}\big)^{\frac{2-d}2}\frac{K_{\frac{d-2}2}(\sqrt{(r-\lambda_0(r,A;\epsilon_0))\frac2D}\thinspace x)}{K_{\frac{d-2}2}(\sqrt{(r-\lambda_0(r,A;\epsilon_0))\frac2D}\thinspace\epsilon_0)}\Big),\ x\ge\epsilon_0.
\end{equation}
\end{2'}
\begin{C1'}
Let
$$
C(r,\epsilon_0,D)=r^{\frac34}\epsilon_0^{\frac{d-2}2}(\frac{\pi^2D}4)^\frac14\big(K_{\frac{d-2}2}(\sqrt{\frac{2r}D}\epsilon_0)\big)^{-1}.
$$
There exist  $C_i(r,A,\epsilon_0,D)$, $i=1,2,3$, satisfying
$$
\lim_{A\to\infty}C_i(r,A,\epsilon_0,D)=1,\ i=1,2,3,
$$
 such that
\begin{equation}\label{2sidedevbddmult}
\begin{aligned}
&C(r,\epsilon_0,D)C_1(r,A,\epsilon_0, D)A^{\frac{1-d}2}
e^{-\sqrt{\frac{2r}D}A}\le\lambda_0(r,A;\epsilon_0)\le\\
& C(r,\epsilon_0,D)C_2(r,A,\epsilon_0,D) A^{\frac{1-d}2} e^{-C_3(r,A,\epsilon_0,D)\sqrt{\frac{2r}D}A}.
\end{aligned}
\end{equation}
\end{C1'}

We now consider  the adjoint semigroup  $\tilde{\mathcal{T}}^{(r,A;\epsilon_0)}_t$  to the semigroup $\mathcal{T}^{(r,A;\epsilon_0)}_t$. Since   $\mathcal{T}^{(r,A;\epsilon_0)}_t$ is  defined on the Banach space
$C_{0_{\epsilon_0}}\big([\epsilon_0,\infty]\big)$,
the adjoint $\tilde{\mathcal{T}}^{(r,A;\epsilon_0)}_t$  operates on the dual space of bounded linear functions on $C_{0_{\epsilon_0}}\big([\epsilon_0,\infty]\big)$.
Since
$[\epsilon_0,\infty]$ with the one-point compactification topology is compact,
this dual space is the space of finite signed measures on $(\epsilon_0,\infty]$ \cite[p.28]{SS}.
 (The reason the measures are on $(\epsilon_0,\infty]$ instead of on $[\epsilon_0,\infty]$ is that $f(\epsilon_0)=0$, for $f\in C_{0_{\epsilon_0}}\big([\epsilon_0,\infty]\big)$.)
Let $\nu$ be such a finite signed measure. We can write $\nu=c_+\nu^+-c_-\nu^-$,
where $\nu^+$ and $\nu^-$ are probability measures on $(\epsilon_0,\infty]$ and $c_+,c_-\ge0$.
From \eqref{Ttmult}, it follows that
\begin{equation}\label{Ttstarmult}
\begin{aligned}
&\tilde{\mathcal{T}}^{(r,A;\epsilon_0)}_t\nu(dy)=c_+\int_{-\infty}^a\nu^+(dx)\mathcal{P}_x^{(r,A)}(Y(t)\in dy;\tau_{\epsilon_0}^{(Y)}>t)-\\
&c_-\int_{-\infty}^a\nu^-(dx)\mathcal{P}_x^{(r,A)}(Y(t)\in dy;\tau_{\epsilon_0}^{(Y)}>t).
\end{aligned}
\end{equation}

Denote the generator of the adjoint semigroup by $\tilde{\mathcal{L}}^{(r,A;\epsilon_0)}$. Of course, this operator has the same principal eigenvalue as does $\mathcal{L}^{(r,A;\epsilon_0)}$.
 \begin{P3'}
The generator   $\tilde{\mathcal{L}}^{(r,A;\epsilon_0)}$ of $\tilde{\mathcal{T}}^{(r,A;\epsilon_0)}_t$
satisfies
\begin{equation}\label{adjointopermult}
 \tilde{\mathcal{L}}^{(r,A;\epsilon_0)}v(y)=\frac D2v''(y)-D\frac{d-1}{2x}v'+D\frac{d-1}{2x^2}v-rv(y)+r\big(\int_{\epsilon_0}^\infty v(x)dx\big)\delta_A(y),\
\end{equation}
for $v$ satisfying $v\in C_{0_{\epsilon_0}}\big([\epsilon_0,\infty]\big)\cap C_b^2((\epsilon_0,\infty))$ and $\int_{\epsilon_0}^\infty|v(x)|dx<\infty$.
Furthermore, a principal eigenfunction $\mathcal{V}_{r,A;\epsilon_0}$  corresponding to the principal eigenvalue $\lambda_0(r,A;\epsilon_0)$ is given by
\begin{equation}\label{adjointpreigenfumult}
\mathcal{V}_{r,A;\epsilon_0}(y)=\begin{cases}
x^{\frac d2}K_{\frac{d-2}2}(qA)
\thinspace\frac{I_{\frac{d-2}2}(q\epsilon_0)K_{\frac{d-2}2}(qx)-
I_{\frac{d-2}2}(qx)K_{\frac{d-2}2}(q\epsilon_0)}
{I_{\frac{d-2}2}(q\epsilon_0)K_{\frac{d-2}2}(qA)-
I_{\frac{d-2}2}(qA)K_{\frac{d-2}2}(q\epsilon_0)},\
 \epsilon_0\le x\le A;\\
%\frac{K_{\frac{d-2}2}}{I_{\frac{d-2}2}}(\sqrt{\frac{2(r-\lambda)}D}A)
x^{\frac d2}K_{\frac{d-2}2}(qx), \ x\ge A,
\end{cases}
% \exp\big(\sqrt{\frac{2(r-\lambda_0(r,a))}D}\thinspace y\big),\ y<0;\\
%\frac{\exp\big(-\sqrt{\frac{2(r-\lambda_0(r,a))}D}\thinspace(y-a)\big)-\exp\big(\sqrt{\frac{2(r-\lambda_0(r,a))}D}\thinspace(y-a)\big)}
%{\exp\big(\sqrt{\frac{2(r-\lambda_0(r,a))}D}\thinspace a\big)-\exp\big(-\sqrt{\frac{2(r-\lambda_0(r,a))}D}\thinspace a\big)}\end{cases},\ 0\le y\le a.
\end{equation}
where $q=\sqrt{\frac{2(r-\lambda_0(r,A;\epsilon_0))}D}$.
 \end{P3'}

\begin{3'}
Let $A>\epsilon_0$. Then
\begin{equation}\label{asympprobmult}
\mathcal{P}_{A}^{(r,A)}(\tau_{\epsilon_0}^{(Y)}>t)=
\frac1{\mathcal{E}_A^{(r,A)}(\mathcal{U}_{r,A;\epsilon_0}(Y(t))|\tau_{\epsilon_0}^{(Y)}>t)}e^{-\lambda_0(r,A;\epsilon_0)\thinspace t},
\end{equation}
where $\mathcal{U}_{r,A;\epsilon_0}$ is as in \eqref{preigenfumult},
Furthermore,
\begin{equation}\label{Mmult}
\begin{aligned}
&\lim_{t\to\infty}\mathcal{E}_A^{(r;A)}(\mathcal{U}_{r,A;\epsilon_0}(Y(t))|\tau_{\epsilon_0}^{(Y)}>t)=
\frac{\int_{\epsilon_0}^\infty \mathcal{U}_{r,A;\epsilon_0}(x)\mathcal{V}_{r,A;\epsilon_0}(x)dx}{\int_{\epsilon_0}^\infty \mathcal{V}_{r,A;\epsilon_0}(x)dx},
\end{aligned}
\end{equation}
where $\mathcal{V}_{r,A;\epsilon_0}$ is as in \eqref{adjointpreigenfumult}.
Thus, for fixed $A$,
\begin{equation}\label{tasympmult}
\mathcal{P}_{A}^{(r;A)}(\tau_{\epsilon_0}^{(Y)}>t)\sim
\frac{\int_{\epsilon_0}^\infty \mathcal{V}_{r,A;\epsilon_0}(x)dx}{\int_{\epsilon_0}^\infty \mathcal{U}_{r,A;\epsilon_0}(x)\mathcal{V}_{r,A;\epsilon_0}(x)dx}
e^{-\lambda_0(r,A; \epsilon_0)\thinspace t},\ \text{as}\ t\to\infty.
\end{equation}
\end{3'}

\begin{P4'}
\begin{equation}\label{atoinfmult}
\lim_{A\to\infty}\mathcal{E}_A^{(r;A)}(\mathcal{U}_{r,A;\epsilon_0}(Y(t))|\tau_{\epsilon_0}^{(Y)}>t)=1,\ \text{uniformly over}\ t\in(0,\infty),
\end{equation}
where $\mathcal{U}_{r,A;\epsilon_0}$ is as in \eqref{preigenfumult}.
\end{P4'}

In the sections that follow, we prove the results stated above in the order that they appeared, except for
Theorem \ref{randomtarget} and Propositions \ref{2dnoreset} and \ref{speed}, whose proofs appear in that order in the final three sections.

\section{Proof of Proposition \ref{Ttcompact}}
We begin by showing that $T_t^{(r,0;a)}$ maps $C_{0_a}([-\infty,a])$ to $C_{0_a}([-\infty,a])$.
Recall that $P_x^{1;(0)}$ and $E_x^{1;(0)}$ denote probabilities and expectations for the  Brownian motion with diffusion parameter $D$ without resetting and started from $x$.
From the definition of the Brownian motion with resetting, we have for $f\in C_{0_a}([-\infty,a])$ and $x\in (-\infty,a)$,
\begin{equation}\label{recursion}
\begin{aligned}
&T^{(r,0;a)}_tf(x)=e^{-rt}E_x^{1;(0)}(f(X(t);\tau_a>t)+\\
&\int_0^tds\thinspace re^{-rs}P_x^{1;(0)}(\tau_a>s)T_{t-s}^{(r,0;a)}f(0).
\end{aligned}
\end{equation}
From this it is easy to see that $T_t^{(r,0;a)}$ maps $C_{0_a}([-\infty,a])$ to $C_{0_a}([-\infty,a])$.
Indeed, it follows readily from standard results that
$\lim_{x\to a}P_x^{1;(0)}(\tau_a>u)=0$, for all $u>0$. From this and \eqref{recursion} it follows that $\lim_{x\to a}T^{(r,0;a)}_tf(x)=0$.
It also follows readily that
for any $N>0$,
$=\lim_{x\to-\infty}P_x^{1;(0)}(X(t)\le -N, \tau_a>t)=1$ and that
 $\lim_{x\to-\infty}P_x^{1;(0)}(\tau_a>s)=1$, for all $s>0$. Using these last two facts, if follows from
\eqref{recursion} that
$\lim_{x\to-\infty}T_t^{(r,0;a)}f(x)$ exists for $f\in C_{0_a}([-\infty,a])$.
Finally, from \eqref{recursion}   it follows that $T_t^{(r,0;a)}f(x)$ inherits its continuity for $x\in(-\infty, a)$
 from the well-known continuity
of  $E_x^{1;(0)}(f(X(t);\tau_a>t)$ and $P_x^{1;(0)}(\tau_a>s)$.
This completes the proof that $T_t^{(r,0;a)}$ maps $C_{0_a}([-\infty,a])$ to $C_{0_a}([-\infty,a])$.

We now show that $T_t^{(r,0;a)}$ is a compact operator.
We write
\begin{equation}\label{termsofdensity}
E_x^{1;(0)}(f(X(t);\tau_a>t)=\int_{-\infty}^ap^{(a)}(t,x,y)f(y)dy,
\end{equation}
where
$p^{(a)}(t,x,y)$ is the transition sub-probability density for the Brownian motion with diffusion parameter $D$ without resetting, and killed upon hitting $a$.
Using the reflection principle, one can show that
\begin{equation}\label{transitiondensity}
\begin{aligned}
&p^{(a)}(t,x,y)=\frac1{\sqrt{2\pi Dt}}\exp(-\frac{(y-x)^2}{2Dt})-\\
&\int_0^t ds\frac{a-x}{\sqrt{2\pi D}s^{\frac32}}
\exp(-\frac{(a-x)^2}{2Ds})
\frac1{\sqrt{2\pi D(t-s)}}\exp(-\frac{(y-a)^2}{2D(t-s)}).
\end{aligned}
\end{equation}
Using \eqref{transitiondensity} along with \eqref{termsofdensity} and \eqref{recursion} shows that
 $T_t^{(r,0;a)}$ maps
bounded sets in $C_{0_a}([-\infty,a])$ to equicontinuous and bounded sets in $C_{0_a}([-\infty,a])$.
This proves the compactness.

We now
turn to the generator.
Let $f\in C^2_b((-\infty,a))\cap\{f:f,L^{1;(r,0)}f\in C_{0_a}\big([-\infty,a]\big)$.
Note that from this assumption, it also follows that $f''\in C([-\infty,a])$.
From \eqref{recursion} and \eqref{termsofdensity}, we have
$$
\begin{aligned}
&\frac1t\big(T_t^{(r,0;a)}f(x)-f(x)\big)=
\frac1te^{-rt}\int_{-\infty}^ap^{(a)}(t,x,y)\big(f(y)-f(x)\big)dy+\\
&\frac1t\int_0^tds\thinspace re^{-rs}P_x^{1;(0)}(\tau_a>s)\big(T_{t-s}^{(r,0;a)}f(0)-f(x)\big).
\end{aligned}
$$
Clearly,
\begin{equation}\label{2ndpart}
\lim_{t\to0}\frac1t\int_0^tds\thinspace re^{-rs}P_x^{(0)}(\tau_a>s)\big(T_{t-s}^{(r,0;a)}f(0)-f(x)\big)=r(f(0)-f(x)).
\end{equation}
Also, from \eqref{transitiondensity}, we have
\begin{equation}\label{1stpart}
\begin{aligned}
&\lim_{t\to0}\frac1te^{-rt}\int_{-\infty}^ap^{(a)}(t,x,y)\big(f(y)-f(x)\big)dy=\\
&\lim_{t\to0}\frac1te^{-rt}\int_{-\infty}^a\frac1{\sqrt{2\pi Dt}}\exp(-\frac{(y-x)^2}{2Dt})\big(f(y)-f(x)\big)dy=\frac D2f''(x).
\end{aligned}
\end{equation}
The first equality in \eqref{1stpart} follows from
the fact that $\int_0^t ds\frac{a-x}{\sqrt{2\pi D}s^{\frac32}}
\exp(-\frac{(a-x)^2}{2Ds})=o(t)$ as $t\to0$.
When the term $e^{-rt}$ is absent, the second equality in \eqref{1stpart}  is the classical calculation
for the generator of  Brownian motion, obtained
by writing $f$ in a Taylor series  with remainder in the form
\begin{equation}\label{taylor}
f(y)=f(x)+f'(x)(y-x)+\frac{f''(x)}2(y-x)^2
+\frac{f''(c_y)-f''(x)}2(y-x)^2.
\end{equation}
 It is easy to show that the equality still holds with $e^{-rt}$ present since this term approaches
 1 when $t\to0$.
From \eqref{2ndpart} and \eqref{1stpart} we obtain
\begin{equation}\label{genconv}
\lim_{t\to0}\frac1t\big(T_t^{(r,0;a)}f(x)-f(x)\big)=\frac D2f''(x)+r(f(0)-f(x))=(L^{1;(r,0)}f)(x).
\end{equation}
By assumption, $L^{1;(r,0)}f\in C_{0_a}([-\infty,a])$.
Furthermore, since $f\in C([-\infty,a])$, it is   uniformly continuous on $[-\infty, a]$, and consequently
it follows   that the convergence in  \eqref{2ndpart} is uniform. Also, since
$f''\in C([-\infty,a])$, it is also uniformly continuous, and thus it follows from \eqref{taylor} that the convergence
with regard to  the second equal sign in \eqref{1stpart} is uniform.
Finally, the fact that $f(a)=0$ guarantees the uniform  convergence to 0 of
the difference between the two expressions on either side of the first equal sign in \eqref{1stpart}.
 Thus, the convergence in \eqref{genconv} is uniform.
 This completes the proof of the calculation of the generator $L^{(r,0;a)}$.
\hfill $\square$

\section{Proofs of Theorem \ref{prinev} and Corollary \ref{cor1}}
\noindent \it Proof of Theorem \ref{prinev}.\rm\
As noted after Proposition \ref{Ttcompact}, $L^{(r,0;a)}$ has a compact resolvent. Thus, by Proposition \ref{Ttcompact} and the Krein-Rutman theorem, it follows that if we find a $\lambda\in\mathbb{R}$ and a function
$u$ satisfying
\begin{equation}\label{evequa}
\begin{aligned}
&\frac D2u''(x)+r(u(0)-u(x))=-\lambda u\ \ \text{in} \ (-\infty,a);\\
&\lim_{x\to-\infty}u(x)\ \text{exists and is finite};\\
&u(a)=0;\\
&u>0 \ \text{in}\ (-\infty,a),
\end{aligned}
\end{equation}
 then $\lambda$ is necessarily the principal eigenvalue $\lambda_0(r,0;a)$, and $u$ is a corresponding principal eigenfunction.
In order to solve  the above nonstandard, homogenous linear equation involving evaluation at a point,
for an appropriate $\lambda$,
 we consider the following
standard, inhomogeneous linear equation involving a free parameter $c\in\mathbb{R}$:
\begin{equation}\label{B}
\begin{aligned}
&\frac D2B_{c,\lambda}''+(\lambda-r)B_{c,\lambda}=-rc,\ x\in(-\infty,a);\\
&\lim_{x\to-\infty}B_{c,\lambda}(x)\ \text{exists and is finite};\\
&B_{c,\lambda}(a)=0;\\
&B_{c,\lambda}>0\ \text{in}\ (-\infty, a).
\end{aligned}
\end{equation}
We will solve explicitly for $B_{c,\lambda}$, for any $c$ and $\lambda$,  and then we look for a solution $(c,\lambda)$ to the equation $B_{c,\lambda}(0)=c$.
Note that if $(c,\lambda)$ solves this equation, then $B_{c,\lambda}$ solves \eqref{evequa}.
It turns out that the set of solutions is of the form $\{(c,\lambda^*): c\in\mathbb{R}\}$, for a unique $\lambda^*$.

Define
\begin{equation}\label{Bform}
\overline{B}_{c,\lambda}=B_{c,\lambda}-\frac{rc}{r-\lambda}.
\end{equation}
Then $B_{c,\lambda}$ solves \eqref{B} if and only if $\overline{B}_{c,\lambda}$ solves
\begin{equation}\label{Bbar}
\begin{aligned}
&\frac D2\overline{B}_{c,\lambda}''+(\lambda-r)\overline{B}_{c,\lambda}=0,\ x\in(-\infty,a);\\
&\lim_{x\to-\infty}\overline{B}_{c,\lambda}(x)\ \text{exists and is finite};\\
&\overline{B}_{c,\lambda}(a)=-\frac{rc}{r-\lambda};\\
&\overline{B}_{c,\lambda}>-\frac{rc}{r-\lambda} \text{in}\ (-\infty, a).
\end{aligned}
\end{equation}
If $\lambda>r$, the general solution to the ODE will involve sines and cosines,  and thus will not satisfy the second line in \eqref{Bbar}.
Thus, we may assume that $\lambda\in(0,r)$.
The general solution to the homogenous ODE in the first line of \eqref{Bbar} is
of the form $C_1\exp(\sqrt{\frac{2(r-\lambda)}D}\thinspace x)+C_2\exp(-\sqrt{\frac{2(r-\lambda)}D}\thinspace x)$.
In light of the requirement in the second line of \eqref{Bbar}, it follows
that  $\overline{B}_{c,\lambda}=C\exp(\sqrt{\frac{2(r-\lambda)}D}\thinspace x)$,
for some $C$. From the third line of \eqref{Bbar}, it follows that
\begin{equation}\label{Bhatsolve}
\overline{B}_{c,\lambda}(x)=-\frac{rc}{r-\lambda}\exp\big(-\sqrt{\frac{2(r-\lambda)}D}\thinspace(a-x)\big).
\end{equation}
Note that  $\overline{B}_{c,\lambda}$  in \eqref{Bhatsolve} also satisfies the fourth line in \eqref{Bbar}.
From \eqref{Bform} and \eqref{Bhatsolve} we obtain
\begin{equation}\label{Bsolve}
B_{c,\lambda}(x)=\frac{rc}{r-\lambda}(1-\exp\big(-\sqrt{\frac{2(r-\lambda)}D}\thinspace(a-x)\big).
\end{equation}

We now solve for $(c,\lambda)=(c^*,\lambda^*)$ in the equation
$B_{c,\lambda}(0)=c$. From \eqref{Bsolve}, this equation gives
\begin{equation}\label{lambdaev}
\lambda =r\exp\big(-a\sqrt{\frac{2(r-\lambda)}D} \big).
\end{equation}
It is easy to check that the  function $\psi(\lambda)=r\exp\big(-a\sqrt{\frac{2(r-\lambda)}D}\big)-\lambda$,
is convex for $\lambda\in[0,r]$. It satisfies $\psi(0)>0$,  $\psi(r)=0$ and $\lim_{\lambda\to r}\psi'(\lambda)=\infty$.
Therefore, there exists a unique $\lambda=\lambda^*\in(0,r)$ that solves \eqref{lambdaev}.
Thus, there exist solutions to \eqref{evequa} if $\lambda=\lambda^*$, and thus
 $\lambda_0(r,0;a)=\lambda^*$. Up to a positive multiplicative constant, the
  solution $u$ to \eqref{evequa} with $\lambda=\lambda_0(r,a)$ is given by \eqref{Bsolve} with $\lambda=\lambda_0(r,0;a)$:
\begin{equation}\label{u}
u(x)=\frac{r}{r-\lambda_0(r,0;a)}(1-\exp\big(-\sqrt{\frac{2\big(r-\lambda_0(r,0;a)\big)}D}\thinspace(a-x)\big),\ x<a.
\end{equation}
This proves \eqref{evequation} and \eqref{preigenfu}. \hfill $\square$

\noindent \it\ Proof of Corollary \ref{cor1}. \rm\
From the fact that $\lambda_0(r,0;a)$ is the unique solution of \eqref{lambdaev} in $(0,r)$,  it follows easily that $\lambda_0(r,0;a)$ is decreasing in $a$.
The corollary follows from this fact along with \eqref{lambdaev}.\hfill $\square$

\section{Proof of Proposition \ref{adjoint}}
By linearity, it suffices to prove \eqref{adjointoper} in the case that $v$, as in the statement of  the theorem, is
a probability density on $(-\infty,a)$; that is, $v\ge0$ and  $\int_{-\infty}^av(x)dx=1$. For such $v$, we need to show that
\begin{equation}\label{adjointgendens}
\tilde{ L}^{(r,0;a)}v(y)=\frac D2v''(y)-rv(y)+r\delta_0(y).
\end{equation}
Recall that
 $P_x^{1;(0)}$ and $E_x^{1;(0)}$ denote probabilities and expectations for the  Brownian motion with diffusion parameter $D$ without resetting and started from $x$.
From \eqref{recursion}, we have
\begin{equation}\label{adjointgencalc}
\begin{aligned}
&\frac1t\big(\tilde T_t^{(r,0;a)}v(y)-v(y)\big)=\frac1te^{-rt}\big(\int_{-\infty}^av(x)p^{(a)}(t,x,y)-v(y)\big)+\\
&\frac1t\int_0^tds\thinspace re^{-rs}\big(P_v^{1;(0)}(\tau_a>s)\tilde T_t^{(r,0;a)}\delta_0(y)-v(y)\big),
\end{aligned}
\end{equation}
where $P_v^{1;(0)}(\tau_a>s)=\int_{-\infty}^av(x)P_x^{1;(0)}(\tau_a>s)dx$.
Clearly,
\begin{equation}\label{secondpart}
\lim_{t\to0}\frac1t\int_0^tds\thinspace re^{-rs}\big(P_v^{1;(0)}(\tau_a>s)\tilde T_t^{(r,0;a)}\delta_0(y)-v(y)\big)=r\big(\delta_0(y)-v(y)\big).
\end{equation}
Also, we have
\begin{equation}\label{firstpart}
\lim_{t\to\infty}\frac1te^{-rt}\big(\int_{-\infty}^av(x)p^{(a)}(t,x,y)-v(y)\big)=\frac D2v''(y),
\end{equation}
by the same argument used for \eqref{1stpart}.
The same argument as at the end of the proof of Proposition \ref{Ttcompact} shows that the convergence in
\eqref{secondpart} and \eqref{firstpart} is uniform.
Thus,
 \eqref{adjointgendens} follows from \eqref{adjointgencalc}-\eqref{firstpart}.

We now turn to obtaining the principal eigenfunction in  \eqref{adjointpreigenfu}.
We need to solve
\begin{equation}\label{needto}
\begin{aligned}
&\frac D2v_0''(y)-(r-\lambda)v_0(y)+r(\int_{-\infty}^av_0(x)dx)\delta_0(y)=0;\\
&v_0(a)=0,\ \ v_0>0 \ \text{in}\ (-\infty,a);\\
& \int_{-\infty}^av_0(y)dy<\infty, \ \lim_{y\to-\infty}v(y)=0,
\end{aligned}
\end{equation}
where $\lambda=\lambda_0(r,0;a)$.
Let  $q=\sqrt{\frac{2(r-\lambda)}D}$. Noting that $e^{qx}$ and $e^{-qx}$ are two linearly independent solutions to the linear ODE obtained from the first line in \eqref{needto} by  deleting the final term on the left hand side involving
the measure $\delta_0$, we look for a solution to \eqref{needto} in the form
\begin{equation}\label{voform}
v_0(y)=\begin{cases} e^{qy},y<0;\\ ce^{-qy}+(1-c)e^{qy},\ 0\le y\le a,\end{cases}
\end{equation}
for some $c\in\mathbb{R}$.
Note that $v_0$ satisfies the ODE in the first line of \eqref{needto} for $y\neq0$.
Also,
$v_0$ is continuous at $y=0$ and  $\int_{-\infty}^av_0(y)dy<\infty$.
In order to obtain $v_0(a)=0$,
we need
\begin{equation}\label{cfromv0form}
c=\frac{e^{qa}}{e^{qa}-e^{-qa}}.
\end{equation}
This completely determines $v_0$ as above, and
plugging $c$ from \eqref{cfromv0form} into \eqref{voform} shows that $v_0\ge0$ and gives
\eqref{adjointpreigenfu}.
However,  we have not yet dealt with the $\delta$-measure in \eqref{needto}.
This is where the particular value $\lambda=\lambda_0(r,0;a)$ comes in. By the Krein-Rutman theory, there must be one (and only one) value of $\lambda$ for which this $v_0$ satisfies
\eqref{needto}.
We could stop here, but
since the calculations are simple, we now verify this explicitly.

We note that for a continuous function $f$ on $(-\infty, a]$ whose second derivative exists  except at $x=0$ and is bounded near $x=0$,   one has
$$
\frac{d}{dy^2}(f(y))=f''(y)+\big(f'(0^+)-f'(0^-)\big)\delta_0(y),
$$
in the sense of distributions. That is,
$$
\int_{-\infty}^a u''(y)f(y)dy=\int_{-\infty}^a u(y)f''(y)dy + (f'(0^+)-f'(0^-))u(0),
$$
for smooth $u$ with compact support in $(-\infty, a)$.
Thus, writing
$v_0(y)=e^{qy}+f(y)$, where
$f(y)=0$, for $y\le0$ and $f(y)=c(e^{-qy}-e^{qy})$, for $y\in[0,a]$,
and noting that $f'(0^+)-f'(0^-)=-2cq$,
we have
\begin{equation}\label{operatorwithdelta}
\big(\frac D2v_0''(y)-(r-\lambda)v_0(y)\big)=-Dcq\delta_0(y)=-\frac{Dqe^{qa}}{e^{qa}-e^{-qa}}\delta_0(y).
\end{equation}
From \eqref{operatorwithdelta}, in order that $v_0$ solve \eqref{needto}, we need
\begin{equation}\label{solveforlambda}
\int_{-\infty}^av_0(y)dy=\frac {Dq}r\frac{e^{qa}}{e^{qa}-e^{-qa}}.
\end{equation}
A direct calculation reveals that
\begin{equation}\label{integralv0}
\int_{-\infty}^av_0(y)dy=\frac2q\frac{e^{qa}-1}{e^{qa}-e^{-qa}}.
\end{equation}
Thus, from \eqref{solveforlambda} and \eqref{integralv0} we need
\begin{equation}\label{Dqrequal}
\frac{Dqe^{qa}}r=\frac2q(e^{qa}-1).
\end{equation}
Recalling that $q=\sqrt{\frac{2(r-\lambda)}D}$, \eqref{Dqrequal} reduces
to $\lambda=re^{-qa}=re^{-a\thinspace \sqrt{\frac{2(r-\lambda)}D}}$.
By \eqref{evequation} in Theorem \ref{prinev}, it follows that \eqref{Dqrequal} holds precisely
for $\lambda=\lambda_0(r,0;a)$.
\hfill $\square$

\section{Proofs of Theorem \ref{Main} and  Proposition \ref{alimitunift}}
\noindent \it Proof of Theorem \ref{Main}.\rm\ We begin with the proof of \eqref{asympprob}. From the standard theory of Markov processes, it follows that
 $$
 f(X(t),t)-\int_0^t(f_t+L^{1;(r,0)}f)(X(s),s)ds
 $$ is a martingale,
 for any $f$ satisfying $f\in C^{2,1}_b((-\infty,\infty)\times(0, T))
 \cap C((-\infty,\infty)\times[0,\infty))$, for all $T>0$,
  where $X(t)$ is the Brownian motion with resetting with generator $L^{1;(r,0)}$ as in \eqref{generator}.
Then by Doob's optional stopping theorem,
$$
f(X(t\wedge \tau_a),t\wedge\tau_a)-\int_0^{t\wedge\tau_a}(f_t+L^{1;(r,0)}f)(X(s),s)ds
$$
 is also a martingale.
Since the process $X(t)$ is stopped at $a$, we can choose $f(x,t)=e^{\lambda_0(r,0;a)t}u_{r,0;a}(x)$, where $u_{r,0;a}$ is as in
\eqref{preigenfu}  and solves \eqref{evequa} with $\lambda=\lambda_0(r,0;a)$. This choice of $f$ gives  $f_t+L^{1;(r,0)}f=0$. Thus,
$$
e^{\lambda_0(r,0;a)(t\wedge\tau_a)}u_{r,0;a}(X(t\wedge\tau_a))\ \text{is a martingale}.
$$
Consequently,
\begin{equation}\label{keyequ}
E_0^{1;(r,0)}e^{\lambda_0(r,0;a)(t\wedge\tau_a)}u_{r,0;a}(X(t\wedge\tau_a))=u_{r,0;a}(0).
\end{equation}
From
\eqref{evequation} and \eqref{preigenfu}, it follows that $u_{r,0;a}(0)=1$. Also,
 $u_{r,0;a}$ vanishes at $a$. Thus, \eqref{keyequ} reduces to
\begin{equation}\label{keyequ1}
e^{\lambda_0(r,0;a)t}E_0^{1;(r,0)}(u_{r,0;a}(X(t);\tau_a>t)=1.
\end{equation}
%Also, since $u$ vanishes at $a$, from now on we can consider $X(t)$ to be the Markov process killed at $a$, with semigroup $T_t^{(r,a)}$ and generator $L^{(r,a)}$.
Writing $E_0^{1;(r,0)}(u_{r,0;a}(X(t);\tau_a>t)=P_0^{1;(r,0)}(\tau_a>t)E_0^{1;(r,0)}(u_{r,0;a}(X(t))|\tau_a>t)$,
we can rewrite \eqref{keyequ1} in the form
\begin{equation}\label{keyequ2}
P_0^{1;(r,0)}(\tau_a>t)=\frac1{E_0^{1;(r,0)}(u_{r,0;a}(X(t))|\tau_a>t)}e^{-\lambda_0(r,0;a)t},
\end{equation}
which is \eqref{asympprob}.

We now turn to the proof of \eqref{M}.
By Proposition \ref{Ttcompact}, the semigroup $T_t^{(r,0;a)}$  is compact. It follows then that
\begin{equation}\label{condconv}
\lim_{t\to\infty}E_0^{1;(r,0)}(u_{r,0;a}(X(t))|\tau_a>t)=
\frac{\int_{-\infty}^au_{r,0;a}(x)v_{r,0;a}(x)dx}{\int_{-\infty}^av_{r,0;a}(x)dx},
\end{equation}
 where $v_{r,0;a}$, appearing  in Proposition \ref{adjoint}, is the principal eigenfunction corresponding to the principal eigenvalue $\lambda_0(r,0;a)$
for
the adjoint operator $\tilde{L}^{(r,0;a)}$.
This follows for example from the corollary after Theorem 3 in \cite{Pin90}.
Actually, that corollary, if it could be applied directly to the situation at hand, would give the stronger result that the transition sub-probability density
$P_x^{1;(r,0)}(X(t)\in dy|\tau_a>t)$ converges uniformly in $x$ and $y$ to
$\frac{v_0(y)}{\int_{-\infty}^av_0(x)dx}$.
However, for the  proof of this as in \cite{Pin90}, we would need to know that
this transition sub-probability density, call it $p^{(r,a)}(t,x,y)$,  satisfies $\sup_{x,y\in(-\infty,a)}p^{(r,a)}(1,x,y)<\infty$. The transition probability in \cite{Pin90} satisfied a standard
parabolic pde, whereas in the situation at hand $p^{(r,a)}(t,x,y)$ satisfies a nonstandard parabolic pde which includes evaluation at 0.
Rather than attempt to prove that the above boundedness condition holds for $p^{(r,a)}(t,x,y)$, we note that in order to prove the weaker form \eqref{condconv}, the method of proof in \cite{Pin90} works
without the necessity of the above uniform pointwise  bound.

Letting $q$ be as in \eqref{M}, and recalling the definition of $v_{r,0;a}$ from \eqref{adjointpreigenfu},  direct calculation reveals that
$$
\int_{-\infty}^av_{r,0;}(x)dx=\frac2q\frac{e^{qa}-1}{e^{qa}-e^{-qa}},
$$
(as has already been noted in \eqref{integralv0}.)
Using \eqref{Dqrequal}, we can rewrite the right hand side above to obtain
\begin{equation}\label{intv0}
\int_{-\infty}^av_{r,0;a}(x)dx=\frac{qDe^{qa}}{r(e^{qa}-e^{-qa})}.
\end{equation}
Recalling also the definition of $u_{r,0;a}$ from \eqref{preigenfu}, a direct calculation gives
\begin{equation}\label{intu0v0}
\begin{aligned}
&\int_{-\infty}^au_{r,0;a}(x)v_{r,0;a}(x)dx=\frac r{r-\lambda_0(r,0;a)}(\frac1q-\frac{e^{-qa}}{2q})+\\
&\frac r{2q(r-\lambda_0(r,0;a))}\frac1{e^{qa}-e^{-qa}}\big(2e^{qa}+2e^{-qa}-e^{-2qa}-3-2qa\big).
\end{aligned}
\end{equation}
After some algebra, \eqref{M} follows from \eqref{intv0} and \eqref{intu0v0}.
Finally, \eqref{tasymp} follows immediately from \eqref{asympprob} and \eqref{M}.
\hfill$\square$
\medskip

\noindent \it Proof of Proposition \ref{alimitunift}.\rm\
Note that for any $y>0$, $u_{r,0;a}$ satisfies
$\lim_{a\to\infty}u_{r,0;a}(x)=1$, uniformly over  $x\in(-\infty, y]$.
Thus, to prove \eqref{atoinf},
it suffices to show that the set of  distributions $\{P_0^{1;(r,0)}(X(t)\in\cdot|\tau_a>t): \ a\ge1, \ t>0\}$
is tight at $+\infty$; namely
\begin{equation}\label{tightinfty}
\lim_{y\to\infty}\sup_{t>0,a\ge1}P_0^{1;(r,0)}(X(t)\ge y|\tau_a>t)=0.
\end{equation}

For each $t>0$, let $LR_t$ be the random variable denoting the last resetting time before time $t$ for the process $X(t)$ under $P_0^{1;(r,0)}$.
Let $\alpha_t(s),\ 0\le s\le t,$ denote the density of the random variable $t-LR_t$, and let $\tilde\alpha_t(s),\ 0\le s\le t,$ denote the density of
$t-LR_t$, when conditioned on $\tau_a>t$.
Recall that $P_0^{1;(0)}$ denotes probabilities for the  Brownian motion without resetting.
From the way the resetting mechanism works, we have
\begin{equation}\label{1stcond}
P_0^{1;(r,0)}(X(t)\ge y|\tau_a>t)=\int_0^t\tilde\alpha_t(s)P_0^{1;(0)}(X(s)\ge y|\tau_a>s)ds.
\end{equation}

We now show that
\begin{equation}\label{2ndcond}
P_0^{1;(0)}(X(s)\ge y|\tau_a>s)\le P_0^{1;(0)}(X(s)\ge y).
\end{equation}
Under $P_0^{1;(0)}$, the process $X(u),0\le u\le s$, conditioned on $\tau_a>s$, is a time-inhomogeneous diffusion process generated
by
%$\frac D2\big(\frac{d^2}{dx^2}+\frac{\frac {\partial P_x^{(0)}(\tau_a>t-s)}{\partial x}}{P_x(\tau_a>t-s)}\big)$
$\frac D2\big(\frac{d^2}{dx^2}+b^{(s)}(u,x)\frac d{dx}\big)$,
where
$b^{(s)}(u,x)=\frac{w_x(s-u,x)}{w(s-u,x)}$ with
$w(u,x)=P_x^{1;(0)}(\tau_a>u),\ x<a$ \cite{Pin85}.
Clearly, $w_x(s-u,x)\le 0$.
Thus, the drift $b^{(s)}$ is non-positive. Now \eqref{2ndcond} follows from this along with the Ikeda-Watanabe comparison theorem \cite{IW}.

From the definitions of $\alpha_t$ and $\tilde\alpha_t$, we have
\begin{equation}\label{alphaalphatilde}
\tilde\alpha_t(s)=\frac{P_0^{1;(0)}(\tau_a>t-s)P_0^{1;(0)}(\tau_a>s)}{P_0^{1;(0)}(\tau_a>t)}\alpha_t(s).
\end{equation}
We have
\begin{equation}\label{st-st}
\frac{P_0^{1;(0)}(\tau_a>t-s)P_0^{1;(0)}(\tau_a>s)}{P_0^{1;(0)}(\tau_a>t)}\le
\frac{P_0^{1;(0)}(\tau_a>\frac t2)}{P_0^{1;(0)}(\tau_a>t)}.
\end{equation}
As is well-known from the  reflection principle,
\begin{equation*}\label{reflectionpr}
P_0^{1;(0)}(\tau_a>u)=2\int_0^{\frac a{\sqrt{Du}}}\frac1{\sqrt{2\pi}}\exp(-\frac{z^2}2)dz,\ a>0.
\end{equation*}
Thus, the right hand side of \eqref{st-st} is bounded in $t$.
Using this with
 \eqref{1stcond}-\eqref{st-st}, we have
\begin{equation}\label{3rdcond}
P_0^{1;(r,0)}(X(t)\ge y|\tau_a>t)\le C\int_0^tP_0^{1;(0)}(X(s)\ge y)\alpha_t(s)ds,
\end{equation}
for some $C>0$.
Thus to prove \eqref{tightinfty},  it suffices to show that
\begin{equation}\label{remainstoshow}
\lim_{y\to\infty}\sup_{t>0}\int_0^\infty\alpha_t(s)P_0^{1;(0)}(X(s)\ge y)ds=0.
\end{equation}
Clearly
\begin{equation}\label{tightnessforBM}
\lim_{y\to\infty}P_0^{1;(0)}(X(s)\ge y)=0, \ \text{uniformly over}\ s\ \text{in a compact set}.
\end{equation}
The distribution with density $\alpha_t$ is obviously stochastically dominated by
the  time that elapses  between the largest resetting time smaller than $t$ and the smallest resetting time larger than $t$.
This distribution is well-known; it has density
$$
f_t(s)=\begin{cases} r^2se^{-rs}, 0\le s\le t;\\ r(1+rt)e^{-rs}, s>t.\end{cases}
$$
\cite[p.13]{F}. Clearly the set of densities
 $\{f_t\}_{0<t<\infty}$ is tight. Thus, the set of  densities $\{\alpha_t\}_{0<t<\infty}$ is tight. From this
and \eqref{tightnessforBM},
it follows that \eqref{remainstoshow} holds.
\hfill $\square$

\section{Proof of Proposition 2-Bessel}\label{2Bessel}
The proof follows similarly to the proof of Proposition \ref{Ttcompact}.
Let $\mathcal{P}_x^{\text{Bes}(d);0}$  and $\mathcal{E}_x^{\text{Bes}(d);0}$ denote probabilities and expectations for the Bessel process of order $d$ and diffusion coefficient $D$ without resetting and starting from $x$.
Similar to \eqref{recursion}, we have
\begin{equation}\label{recursionmult}
\begin{aligned}
&\mathcal{T}^{(r,A;\epsilon_0)}_tf(x)=e^{-rt}\mathcal{E}_x^{\text{Bes}(d);0}(f(Y(t);\tau^{(Y)}_{\epsilon_0}>t)+\\
&\int_0^tds\thinspace re^{-rs}\mathcal{P}_x^{\text{Bes}(d);0}(\tau^{(Y)}_{\epsilon_0}>s)\mathcal{T}_{t-s}^{(r,A;\epsilon_0)}f(0).
\end{aligned}
\end{equation}
Using \eqref{recursionmult}, the proof that
the semigroup operator $\mathcal{T}^{(r,A;\epsilon_0)}_t$ maps
$C_{0_{\epsilon_0}}\big([\epsilon_0,\infty]\big)$ to $C_{0_{\epsilon_0}}\big([\epsilon_0,\infty]\big)$
is just like the  corresponding proof in Proposition \ref{Ttcompact}, because the basic properties of the Brownian motion semigroup that were used in the proof are shared by the Bessel process semigroup.
This is also true with regard to the calculation of the generator.

With regard to the proof that the operator is compact, we rewrite \eqref{recursionmult} as
\begin{equation}\label{recursionmultden}
\begin{aligned}
&\mathcal{T}^{(r,A;\epsilon_0)}_tf(x)=e^{-rt}\int_{\epsilon_0}^\infty p^{\text{Bes}(d);\epsilon_0}(t,x,y)f(y)dy+\\
&\int_0^tds\thinspace re^{-rs}\mathcal{P}_x^{\text{Bes}(d);0}(\tau^{(Y)}_{\epsilon_0}>s)\mathcal{T}_{t-s}^{(r,A;\epsilon_0)}f(0),
\end{aligned}
\end{equation}
where $p^{\text{Bes}(d);\epsilon_0}(t,x,y)$ denotes
 the  transition sub-probability density for the Bessel   process without resetting and  killed upon hitting $\epsilon_0$.
From \eqref{recursionmultden}, it suffices to show that the functions
$\{p^{\text{Bes}(d);\epsilon_0}(t,\cdot,y)\}_{y\in(\epsilon_0,\infty)}$ and
the functions $\{\mathcal{P}_\cdot^{\text{Bes}(d);0}(\tau^{(Y)}_{\epsilon_0}>s)\}_{s>0}$ are uniformly
equicontinuous on $(\epsilon_0,\infty)$. We sketch how
this equicontinuity can be deduced from \cite{BM} and \cite{BMR}.

In those two papers, the parameter  $\mu$  plays the role of    our $\frac d2-1$ (or equivalently, $2\mu+2$ plays the role of our $d$), and in \cite{BM} the parameter $a$ plays the role of our $\epsilon_0$. Also, our $D$ is equal to 1 in those papers. So we only consider this case.
(The general case follows by scaling.)
In \cite{BM},  $\{p_a^\mu(t,\cdot,y)\}_{y\in(a,\infty)}$ plays the role of our $\{p^{\text{Bes}(d),\epsilon_0}(t,\cdot,y)\}_{y\in(\epsilon_0,\infty)}$
and $\{q_{\cdot,a}^\mu(s)\}_{s>0}\}$ plays the role of our  $\{\mathcal{P}_\cdot^{\text{Bes}(d);0}(\tau^{(Y)}_{\epsilon_0}>s)\}_{s>0}$.
So we need to demonstrate the uniform equicontinuity of
$\{p_a^\mu(t,\cdot,y)\}_{y\in(a,\infty)}$ and of $\{q_{\cdot,a}^\mu(s)\}_{s>0}\}$ over $(\epsilon_0,\infty)$.

From  \cite[(2.10)]{BM},
  the uniform equicontinuity of
$\{p_a^\mu(t,\cdot,y)\}_{y\in(a,\infty)}$ follows from that of $\{q_{\cdot,a}^\mu(s)\}_{s>0}\}$ and
of $\{p^\mu(t,\cdot,y)\}_{y\in(a,\infty)}$, where $p^\mu(t,x,y)$ is the transition probability function for the Bessel process
of order $d=2\mu+2$ without killing.
Using \cite[equations (1.1) and (2.3)]{BM} shows the uniform equicontinuity of $\{p^\mu(t,\cdot,y)\}_{y\in(a,\infty)}$.
The uniform equicontinuity of $\{q_{\cdot,a}^\mu(s)\}_{s>0}\}$  can be deduced from section 2.3 in \cite{BMR}.
\hfill $\square$

\section{Proof of Theorem 2-Bessel and Corollary 1-Bessel}\label{Thm1-Bessel}
\noindent \it Proof of Theorem 2-Bessel.\rm\ The proof follows the contours of the  proof of Theorem \ref{prinev}, except that instead of using the operator $\frac D2\frac{d^2}{dx^2}$ on $(-\infty,a]$, we use the operator
$\frac D2\frac{d^2}{dx^2}+\frac{D(d-1)}{2x}\frac d{dx}$ on $[\epsilon_0,\infty)$, and instead of the distinguished point being 0, it is $A$.
We use the same notation $u$, $B_{c,\lambda}$ and $\overline{B}_{c,\lambda}$ for the functions appearing in that proof.
As noted after Proposition 2-Bessel, $\mathcal{L}^{(r,A;\epsilon_0)}$ has a compact resolvent. Thus, by Proposition 2-Bessel and the Krein-Rutman theorem, it follows that the principal eigenvalue
$\lambda_0(r,A;\epsilon_0)$ is the unique solution $\lambda$ to the
 following equation, analogous to  \eqref{evequa}:
\begin{equation}\label{evequamult}
\begin{aligned}
&\frac D2u''(x)+D\frac{(d-1)}{2x}u'(x)+r(u(A)-u(x))=-\lambda u\ \ \text{in} \ (\epsilon_0,\infty);\\
&\lim_{x\to\infty}u(x)\ \text{exists and is finite};\\
&u(\epsilon_0)=0;\\
&u>0 \ \text{in}\ (\epsilon_0,\infty).
\end{aligned}
\end{equation}
And
the function $\overline{B}_{c,\lambda}$ satisfies the following equation, analogous to
\eqref{Bbar}:
\begin{equation}\label{Bbarmult}
\begin{aligned}
&\frac D2\overline{B}_{c,\lambda}''+D\frac{(d-1)}{2x}\overline{B}_{c,\lambda}'+(\lambda-r)\overline{B}_{c,\lambda}=0,\ x\in(\epsilon_0,\infty);\\
&\lim_{x\to\infty}\overline{B}_{c,\lambda}(x)\ \text{exists and is finite};\\
&\overline{B}_{c,\lambda}(\epsilon_0)=-\frac{rc}{r-\lambda};\\
&\overline{B}_{c,\lambda}>-\frac{rc}{r-\lambda}\ \text{in}\ (\epsilon_0,\infty).
\end{aligned}
\end{equation}
The modified Bessel functions of the first and second kind, of order $\nu$, denoted respectively by $I_\nu$ and $K_\nu$ are linearly independent solutions to the linear ODE
$x^2\frac{d^2 W}{dx^2}+x\frac{dW}{dx}-(x^2+\nu^2)W=0$ \cite{AS,W}. By looking for solutions of the form $x^\gamma K_\nu(\eta x)$
and $x^\gamma I_\nu(\eta x)$,
 for parameters $\gamma,\nu$ and $\eta$,
 one can verify that two linearly independent solutions to the ODE in \eqref{Bbarmult} are:
$$
x^{\frac{2-d}2}K_\frac{d-2}2(\sqrt{\frac2D(r-\lambda)}\thinspace x),\ \ \ \
x^{\frac{2-d}2}I_\frac{d-2}2(\sqrt{\frac2D(r-\lambda)}\thinspace x).
$$
The function $I_\nu$ grows exponentially and the function $K_\nu$ decays exponentially as $x\to\infty$ \cite{AS,W}.
Therefore, it follows from \eqref{Bbarmult} that
$$
\overline{B}_{c,\lambda}=-\frac{rc}{r-\lambda}\thinspace\frac{\epsilon_0^\frac{d-2}2}{K_{\frac{d-2}2}(\sqrt{\frac2D(r-\lambda)}\thinspace \epsilon_0)}x^{\frac{2-d}2}K_\frac{d-2}2(\sqrt{\frac2D(r-\lambda)}\thinspace x).
$$
Similar to  the passage from \eqref{Bhatsolve} to \eqref{Bsolve}, we have
\begin{equation}\label{Bsolvemult}
B_{c,\lambda}(x)=\frac{rc}{r-\lambda}\Big(1-\big(\frac {x}{\epsilon_0}\big)^{\frac{2-d}2}\frac{K_{\frac{d-2}2}(\sqrt{(r-\lambda)\frac2D}\thinspace x)}{K_{\frac{d-2}2}(\sqrt{(r-\lambda)\frac2D}\thinspace\epsilon_0)}\Big)
\end{equation}

Similar to   the proof of Theorem \ref{prinev}, we
 now solve for $(c,\lambda)=(c^*,\lambda^*)$ in the equation
$B_{c,\lambda}(A)=c$. From \eqref{Bsolvemult}, this yields the equation
\eqref{evequationmult} for $\lambda=\lambda^*$, with
$c=c^*$ being arbitrary.
Although we could perform an analysis to show directly that there exists a unique $\lambda^*\in(0,r)$ that solves \eqref{evequationmult},
similar to what was done in the proof of Theorem \ref{prinev}, but much more tedious, in the present case we simply note that this follows by the uniqueness of the principal eigenvalue in
the Krein-Rutman theorem. This unique solution is the principal eigenvalue $\lambda_0(r,A;\epsilon_0)$.
A principal eigenfunction $\mathcal{U}_{0;A}(x)$ is then given by the right hand side of \eqref{Bsolvemult}
with $\lambda=\lambda_0(r,A;\epsilon_0)$ and, say, $c=1$.
\hfill $\square$
\medskip

\noindent \it Proof of Corollary 1-Bessel.\rm\ The corollary follows readily from \eqref{evequationmult} and the asymptotic estimate
$K_{\frac{d-2}2}(A)\sim\sqrt{\frac\pi{2A}}e^{-A}$ as $A\to\infty$ \cite{AS,W}.
(Note that this leading order asymptotic behavior for  $K_{\frac{d-2}2}(A)$ is independent of the order $\frac{d-2}2$.)
\hfill $\square$

\section{Proof of Proposition 3-Bessel}\label{Thm3-Bessel}
The proof that the adjoint generator is as in \eqref{adjointopermult} is similar to the proof of  Proposition \ref{adjoint}, so we leave it to the reader. We turn to the calculation of the corresponding principal eigenfunction
in \eqref{adjointpreigenfumult}.
We need to solve
\begin{equation}\label{needtomult}
\begin{aligned}
&\frac D2v''(y)-D\frac{d-1}{2x}v'+D\frac{d-1}{2x^2}v-(r-\lambda)v(y)+r\big(\int_{\epsilon_0}^\infty v(x)dx\big)\delta_A(y)=0;\\
&v(\epsilon_0)=0, v>0\ \text{in}\ (\epsilon_0,\infty);\\
&\int_{\epsilon_0}^\infty v(y)dy<\infty,\  \lim_{y\to\infty}v(y)=0,
\end{aligned}
\end{equation}
where $\lambda=\lambda(r,A;\epsilon_0)$.
From the proof of Theorem 2-Bessel, recall $K_\nu$ and $I_\nu$, the modified Bessel functions of order $\nu$, which are linearly independent solutions to the linear ODE
$x^2\frac{d^2 W}{dx^2}+x\frac{dW}{dx}-(x^2+\nu^2)W=0$.
%Consider the eigenvalue equation for the  operator $\mathcal{O}$, obtained from $\tilde{\mathcal{L}}^{(r,\epsilon_0,A)}$ by ignoring its final term in which appears $\delta_A$:
%$$
%\mathcal{O}v:=\frac D2v''(y)-D\frac{d-1}{2x}v'+D\frac{d-1}{2x^2}v-rv(y)=-\lambda v.
%$$
By looking for solutions of the  form $x^\gamma K_\nu(\eta x)$, for parameters $\gamma,\nu$ and $\eta$,
one can verify that two linearly independent solutions to the  linear ODE obtained from the first line of \eqref{needtomult} by deleting the final term
involving the  measure $\delta_A$ are
$$
x^{\frac d2}K_{\frac{d-2}2}(\sqrt{\frac2D(r-\lambda)}\thinspace x),\ \ \ \ \ \ x^{\frac d2}I_{\frac{d-2}2}(\sqrt{\frac2D(r-\lambda)}\thinspace x).
$$
Recalling that $I_{\frac{d-2}2}$ grows exponentially and $K_{\frac{d-2}2}$ decays exponentially,
we look for the solution $v$ to \eqref{needtomult} in the form
\begin{equation}\label{formofsolu}
\begin{aligned}
&v(x)=\begin{cases}c_1x^{\frac d2}I_{\frac{d-2}2}(\sqrt{\frac2D(r-\lambda)}\thinspace x)+c_2x^\frac d2K_{\frac{d-2}2}(\sqrt{\frac2D(r-\lambda)}\thinspace x),\\
 \epsilon_0\le x\le A;\\
x^\frac d2K_{\frac{d-2}2}(\sqrt{\frac2D(r-\lambda)}\thinspace x),\ x\ge A.\end{cases}
\end{aligned}
\end{equation}
Then $v$ satisfies the third line of \eqref{needtomult} and it satisfies the ODE in the first line of \eqref{needtomult} for $y\neq A$.
Using the equations $v(\epsilon_0)=0$ and $v(A^-)=v(A^+)$, we can solve for $c_1$ and $c_2$, obtaining
 \eqref{adjointpreigenfumult}.
The $\delta$-measure requirement at $y=A$ in the first line of \eqref{needtomult} follows automatically from the Krein-Rutman theorem in the case that $\lambda=\lambda(r,A;\epsilon_0)$.
(See the discussion at the corresponding juncture of the proof of Proposition \ref{adjoint}, which contains the corresponding result in the one-dimensional case.)

\section{Proofs of Theorem 3-Bessel and Proposition 4-Bessel}\label{thm3-Bessel}
\noindent \it Proof of Theorem 3-Bessel.\rm\ The proof is just like the proof of Theorem
\ref{Main}.\hfill \ $\square$

\noindent \it Proof of Proposition 4-Bessel.\rm\ From \eqref{preigenfumult} we have
$\lim_{x,A\to\infty}\mathcal{U}_{r,A;\epsilon}(x)=1$.  Thus, to prove \eqref{atoinfmult},
 it suffices to prove that
\begin{equation}\label{sufficesmathcalU}
\lim_{y\to\infty}\limsup_{A\to\infty}\sup_{t>0}\mathcal{P}_A^{(r;A)}(Y(t)\le A-y)|\tau_{\epsilon_0}^{(Y)}>t)=0.
\end{equation}
For each $t>0$, let
$\text{LR}_t$ be the random variable denoting the last resetting time before $t$ for the process $Y(\cdot)$ under $\mathcal{P}_A^{(r;A)}$.
Let $\alpha_t(s), 0\le s\le t$, denote the density of the random variable $t-\text{LR}_t$, and let $\tilde\alpha_t(s), 0\le s\le t$, denote the density of
$t-\text{LR}_t$ when conditioned on $\tau_{\epsilon_0}^{(Y)}>t$. From the way the resetting mechanism works, we have
\begin{equation}\label{firstcondmult}
\mathcal{P}_A^{(r;A)}(Y(t)\le A-y)|\tau_{\epsilon_0}^{(Y)}>t)=\int_0^t\tilde\alpha_t(s)\mathcal{P}_A^{(r;A)}(Y(s)\le A-y|\tau_{\epsilon_0}^{(Y)}>s)ds.
\end{equation}

We now show that
\begin{equation}\label{2ndcondmult}
\mathcal{P}_A^{(r;A)}(Y(s)\le A-y|\tau_{\epsilon_0}^{(Y)}>s)\le\mathcal{P}_A^{(r;A)}(Y(s)\le A-y).
\end{equation}
Under $\mathcal{P}_A^{(r;A)}$, the process $Y(u), 0\le u\le s$, conditioned on  $\tau_{\epsilon_0}^{(Y)}>s$, is a time-inhomogeneous diffusion process generated by
$\frac D2(\frac{d^2}{dy^2}+\frac{d-1}y\frac d{dy}+ b^{(s)}(u,y)\frac d{dy})$, where
$b^{(s)}(u,y)=\frac{w_y(s-u,y)}{w(s-u,y)}$ with $w(u,y)=\mathcal{P}_y^{(r;A)}(\tau_{\epsilon_0}^{(Y)}>u)$ \cite{Pin85}.
Clearly, $w_y(s-u,y)\ge0$.
Thus, the drift $b^{(s)}$ is nonnegative. Now \eqref{2ndcondmult} follows from this along with the Ikeda-Watanabe comparison theorem \cite{IW}.

From the definitions of $\alpha_t$ and $\tilde\alpha_t$, we have
\begin{equation}\label{alphaalphatildemult}
\tilde\alpha_t(s)=\frac{\mathcal{P}_A^{(r;A)}(\tau_{\epsilon_0}^{(Y)}>t-s))\mathcal{P}_A^{(r;A)}(\tau_{\epsilon_0}^{(Y)}>s)}{\mathcal{P}_A^{(r;A)}(\tau_{\epsilon_0}^{(Y)}>t)}\alpha_t(s).
\end{equation}
We have
\begin{equation}\label{st-stmult}
\frac{\mathcal{P}_A^{(r;A)}(\tau_{\epsilon_0}^{(Y)}>t-s))\mathcal{P}_A^{(r;A)}(\tau_{\epsilon_0}^{(Y)}>s)}{\mathcal{P}_A^{(r;A)}(\tau_{\epsilon_0}^{(Y)}>t)}\le
\frac{\mathcal{P}_A^{(r;A)}(\tau_{\epsilon_0}^{(Y)}>\frac t2)}{\mathcal{P}_A^{(r;A)}(\tau_{\epsilon_0}^{(Y)}>t)}.
\end{equation}
For $d\ge3$, the Bessel process of order $d$ is transient,  so $\lim_{t\to\infty}\mathcal{P}_A^{(r;A)}(\tau_{\epsilon_0}^{(Y)}>t)>0$,
and thus the right hand side of \eqref{st-stmult} is bounded in $t$. For $d=2$, $\mathcal{P}_A^{(r;A)}(\tau_{\epsilon_0}^{(Y)}>t)$ has logarithmic decay (\cite[p.224]{Pin95}) from which it follows that the right hand side
of \eqref{st-stmult} is bounded in $t$. Using this with \eqref{firstcondmult}-\eqref{st-stmult}, we have
\begin{equation}\label{3rdcondmult}
\mathcal{P}_A^{(r;A)}(Y(t)\le A-y)|\tau_{\epsilon_0}^{(Y)}>t)\le C\int_0^t\mathcal{P}_A^{(r;A)}(Y(s)\le A-y)\alpha_t(s)ds,
\end{equation}
for some $C>0$.
Since $Y(\cdot)$ is a Bessel process of order $d$, it is clear that
\begin{equation}\label{additional}
\lim_{y\to\infty}\limsup_{A\to\infty}\mathcal{P}_A^{(r;A)}(Y(s)\le A-y)=0, \ \text{for all} \ s>0.
\end{equation}
As shown at the end of the proof of Proposition \ref{alimitunift},
the distributions $\{\alpha_t\}_{0\le t<\infty}$ are tight.
Now   \eqref{sufficesmathcalU} follows from this along with \eqref{3rdcondmult} and \eqref{additional}.
\hfill $\square$

\section{Proof of Theorem \ref{randomtarget}}
We first prove the one-dimensional case.
In light of \eqref{targetdist}, it suffices to prove \eqref{asympprobtarget} with the range of integration from $0$ to $\infty$ instead of from $-\infty$ to $\infty$.
From \eqref{targetdist} and
\eqref{asympprob}, we have
\begin{equation}\label{intasymp}
\begin{aligned}
&\int_0^\infty  P_0^{1;(r,0)}(\tau_a>t)\mu_{B,l}(da)=\int_0^\infty \frac{e^{-\lambda_0(r,0;a)\thinspace t}}{E_0^{1;(r,0)}(u_{r,0;a}(X(t))|\tau_a>t)}c(a)e^{-Ba^l}da,\\
& \text{as}\ t\to\infty,
\end{aligned}
\end{equation}
where
\begin{equation}\label{casymp}
\lim_{a\to\infty}\frac{\log c(a)}{a^l}=0.
\end{equation}
Since $\lambda_0(r,0;a)$ is decreasing to 0 as $a\to\infty$, it is clear that the asymptotic behavior of the right hand side of \eqref{intasymp} as $t\to\infty$ depends only on large $a$.
Thus, in light of
\eqref{atoinf},
\begin{equation}\label{intasympagain}
\int_0^\infty \frac{e^{-\lambda_0(r,0;a)\thinspace t}}{E_0^{1;(r,0)}(u_{r,0;a}(X(t))|\tau_a>t)}c(a)e^{-Ba^l}da\sim
\int_0^\infty c(a)e^{-\lambda_0(r,0;a)\thinspace t-Ba^l}da,\ \text{as}\ t\to\infty.
\end{equation}
By \eqref{evequation}, we can replace $\lambda_0(r,0;a)$ in the exponent on the right hand side of \eqref{intasympagain} by
$re^{-\sqrt{\frac{2(r-\lambda_0(r,0;a))}D}\thinspace a}$. Making this replacement, using the fact that $\lambda(r,0;a)$ approaches 0 as $a\to\infty$, and using
 \eqref{casymp}, it follows  that if
%\begin{equation}\label{anotherreduction}
%\log \int_0^\infty c(a)e^{-\lambda_0(r,0;a)\thinspace t-Ba^l}da \sim\log\int_0^\infty %\exp(-rte^{-\sqrt{\frac{2r}D}\thinspace a}-Ba^l)da, \ \text{as}\ t\to\infty.
%\end{equation}
\begin{equation}\label{needed}
\lim_{t\to\infty}\frac1{(\log t)^l}\log  \int_0^\infty \exp(-Rte^{-\kappa a}-Ba^l)da=-\frac B{\kappa^l},
\ \text{for all}\ B,R,\kappa>0,
\end{equation}
then
\begin{equation}\label{wantit}
\lim_{t\to\infty}\frac1{(\log t)^l}\log  \int_0^\infty c(a)e^{-\lambda_0(r,0;a)t-Ba^l}da=-B(\frac D{2r})^\frac l2.
\end{equation}
(In fact, it is unnecessary here to replace the specific $r$ with the generic $R$; however, we will need  this general form of \eqref{needed}  in
the proof of the multi-dimensional case.)
Therefore, from   \eqref{intasymp}, \eqref{intasympagain} and \eqref{wantit},  it follows
that the proof of \eqref{asympprobtarget}
will be completed if we prove \eqref{needed}.

To analyze the left hand side of \eqref{needed}, we  locate, for each large $t$, the minimum
 of the expression
\begin{equation}\label{gammadef}
\gamma_t(a):=Rte^{-\kappa a}+Ba^l.
\end{equation}
First consider the case that $l\ge1$.  In this case,
$\gamma_t$ is convex and $\gamma_t'(0)<0$, for all sufficiently large $t$ (actually all $t$, if $l>1$).  Thus, for large $t$, it has a unique minimum which occurs at some $a^*$
which satisfies
\begin{equation}\label{a*}
\kappa Rte^{-\kappa a^*}=lB(a^*)^{l-1}.
\end{equation}
Substituting from \eqref{a*}, we have
\begin{equation}\label{gammaa*}
\gamma_t(a^*)=Rte^{-\kappa\thinspace a^*}+B(a^*)^l=B(a^*)^{l-1}\big(\frac l\kappa+a^*\big).
\end{equation}
From \eqref{a*} we have
$$
\log(\kappa R)+\log t-\kappa a^*=\log(lB)+(l-1)\log a*,
$$
from which it follows that
\begin{equation}\label{a*tasymp}
a^*\sim\frac1\kappa\log t,\ \text{as}\ t\to\infty.
\end{equation}
Substituting from \eqref{a*tasymp} into the right hand side of \eqref{gammaa*},  we have
\begin{equation}\label{finalge1}
\gamma_t(a^*)\sim B(a^*)^l=B(\frac{\log t}\kappa)^l, \ \text{as}\ t\to\infty.
\end{equation}

Now consider the case $l\in(0,1)$.
We have $\gamma_t'(a)=-\kappa Rte^{-\kappa a}+lBa^{l-1}$.
Note that $\gamma_t'(0^+)=\infty$, and it is easy to see that for each $t$, $\gamma_t'(a)>0$ for sufficiently large $a$.
However, $\gamma_t'(1)<0$, for sufficiently large $t$.
Thus, for sufficiently large $t$, there must be at least two roots to $\gamma_t'(a)=0$.
Substituting $\kappa Rte^{-\kappa a}=lBa^{l-a}$ into $\gamma_t''(a)$,
it follows that if $\gamma_t'(a)=0$, then $\gamma_t''(a)=Bla^{l-2}\big(\kappa a+l-1\big)$. Thus, a zero $a$ of $\gamma_t'$ is
a relative maximum of $\gamma_t$ if $a<\frac{1-l}\kappa$ and is a relative minimum if $a>\frac{1-l}\kappa$.
From this it follows that for sufficiently large $t$, there are exactly two zeroes of $\gamma_t'$, and that the larger one is the global minimum of $\gamma_t$,
Denote this global minimum by $a^*$. Using the fact that $a^*>\frac{1-l}\kappa$ and that
$0=\gamma_t'(a^*)=-\kappa Rte^{-\kappa a^*}+lB(a^*)^{l-1}$, it follows that $a^*$ approaches $\infty$ as $t\to\infty$.
The rest of the analysis is as before.

Thus, for sufficiently large $t$,
\eqref{a*tasymp} and \eqref{finalge1} hold for all $l>0$, and
 from the previous  paragraph, for all $l>0$ we have
\begin{equation}\label{2ndderiv}
\gamma_t''(a^*)=Bl(a^*)^{l-2}\big(\kappa a^*+l-1\big).
\end{equation}
Note that
\begin{equation}\label{gamma2ndderiv}
\gamma_t''(a)=\kappa^2Rte^{-\kappa a}+l(l-1)Ba^{l-2}.
\end{equation}
%Now let $a\in(a^*, a^*+1)$ and assume that $t$ is sufficiently large so that $a^*\ge1$.
Using \eqref{2ndderiv} and \eqref{gamma2ndderiv}, we have
\begin{equation}\label{2ndderivcalc}
\begin{aligned}
&\gamma_t''(a)=\gamma_t''(a^*)+(\gamma_t''(a)-\gamma_t''(a^*))=Bl(a^*)^{l-2}\big(\kappa a^*+l-1\big)+\\
&\kappa^2Rt(e^{-\kappa a}-e^{-\kappa a^*})+l(l-1)B(a^{l-2}-(a^*)^{l-2})\le\\
&Bl(a^*)^{l-2}\big(\kappa a^*+l-1\big)+l(l-1)B(a^{l-2}-(a^*)^{l-2}),\ \text{for}\ a>a^*.
\end{aligned}
\end{equation}
From \eqref{2ndderivcalc} and \eqref{a*tasymp} it follows that for some constant $C>0$,
\begin{equation}\label{2ndderivbound}
\gamma''(a)\le C\log^{l-1}t,\  \text{for}\ a\in[a^*,a^*+1].
\end{equation}
Since $\gamma_t'(a^*)=0$, it follows from \eqref{2ndderivbound} that
\begin{equation}\label{gammatbound}
\gamma_t(a)\le\gamma_t(a^*)+\frac12C(\log^{l-1}t)(a-a^*)^2,\ \text{for}\ a\in[a^*,a^*+1].
\end{equation}
Thus, we conclude from \eqref{gammatbound} that for some $\alpha>0$,
\begin{equation}\label{gammatboundagain}
\gamma_t(a)\le\gamma_t(a^*)+1,\ \text{for}\ \begin{cases}\ a\in[a^*,a^*+\alpha], \ \text{if}\ l\in(0,1];\\
a\in\big[a^*,a^*+\frac \alpha{\log^{l-1}t}],\ \text{if}\ l>1.\end{cases}
\end{equation}
Note that the two cases of the interval appearing on the right hand side of \eqref{gammatboundagain} can be merged by writing
$\big[a^*,a^*+\frac{\alpha}{(\log t)^{\max(0,l-1)}}]$.
From this observation along with
 \eqref{gammatboundagain}, \eqref{finalge1} and the definition of $\gamma_t$ in \eqref{gammadef},
 we obtain the lower bound
\begin{equation}\label{lowerfinal}
\begin{aligned}
& \int_0^\infty \exp(-Rte^{-\kappa a}-Ba^l)da\ge \frac\alpha{(\log t)^{\max(0,l-1)}}\exp\big(-(1+\epsilon)B(\frac{\log t}\kappa)^l-1\big),\\
& \ \text{for any}\ \epsilon>0\ \text{and for sufficiently large}\ t\ \text{depending on}\ \epsilon.
\end{aligned}
\end{equation}

Now we turn to an upper bound   for the left hand side of \eqref{lowerfinal}.
Applying L'H\^opital's rule to $\frac{\int_x^\infty e^{-Ba^l}da}{x^{-l+1}e^{-Bx^l}}$ shows
that
\begin{equation*}
\int_x^\infty e^{-Ba^l}da\sim\frac{x^{-l+1}e^{-Bx^l}}{lB},\ \text{as}\ x\to\infty,
\end{equation*}
and thus,
\begin{equation}\label{hopital}
\int_x^\infty e^{-Ba^l}da\le \frac{x^{-l+1}e^{-(1-\epsilon)Bx^l}}{lB},\ \text{for any}\ \epsilon>0\ \text{and sufficiently large}\ x\ \text{depending on}\ \epsilon.
\end{equation}
Write
\begin{equation}\label{inttwoparts}
\begin{aligned}
&\int_0^\infty \exp(-Rte^{-\kappa a}-Ba^l)da=\\
&\int_0^{a^*}\exp(-Rte^{-\kappa a}-Ba^l)da+\int_{a^*}^\infty \exp(-Rte^{-\kappa a}-Ba^l)da.
\end{aligned}
\end{equation}
Using \eqref{hopital} and \eqref{a*tasymp}  gives
\begin{equation}\label{tailpiece}
\begin{aligned}
&\int_{a^*}^\infty \exp(-Rte^{-\kappa a}-Ba^l)da\le
\int_{a^*}^\infty e^{-Ba^l}da\le
 \frac1{lB}\big(\frac{\log t}\kappa\big)^{-l+1}e^{-(1-\epsilon)B\big(\frac{\log t}\kappa\big)^l},\\
&\text{for any}\ \epsilon>0\ \text{and sufficiently large}\ t\ \text{depending on}\ \epsilon.
\end{aligned}
\end{equation}
From the definition of $\gamma_t$ in  \eqref{gammadef} and     the fact that $a^*$ is the minimum of $\gamma_t(a)$,  it follows from \eqref{a*tasymp} and \eqref{finalge1}
that
\begin{equation}\label{firstpiece}
\begin{aligned}
&\int_0^{a^*}\exp(-Rte^{-\kappa  a}-Ba^l)da\le (1+\epsilon)\frac{\log t}\kappa e^{-(1-\epsilon)B(\frac{\log t}\kappa)^l},\\
&\text{for any}\ \epsilon>0\ \text{and sufficiently large}\ t\ \text{depending on}\ \epsilon.
\end{aligned}
\end{equation}
From \eqref{inttwoparts}-\eqref{firstpiece}, we conclude that
\begin{equation}\label{upperfinal}
\begin{aligned}
&\int_0^\infty \exp(-Rte^{-\kappa a}-Ba^l)da\le\\
&(1+\epsilon)\frac{\log t}\kappa e^{-(1-\epsilon)B(\frac{\log t}\kappa)^l}+ \frac1{lB}\big(\frac{\log t}\kappa\big)^{-l+1}e^{-(1-\epsilon)B\big(\frac{\log t}\kappa\big)^l},\\
&\text{for any}\ \epsilon>0\ \text{and sufficiently large}\ t\ \text{depending on}\ \epsilon.
\end{aligned}
\end{equation}
Now \eqref{needed} follows from \eqref{lowerfinal} and \eqref{upperfinal}.
This completes the proof of the one-dimensional case.

We now turn to the multi-dimensional case, where we will also utilize \eqref{needed}.
For $A>0$, let $\bar\mu_{B,l}^{(d)}(A)=\int_{|x|=1}\mu_{B,l}^{(d)}(Ax)s_d(dx)$,
where $s_d$ denotes Lebesgue measure on the unit sphere in $R^d$. We
 have
\begin{equation}\label{intasympBess}
\begin{aligned}
&\int_{\mathbb{R}^d}  P_0^{d;(r,0)}(\tau_a>t)\mu_{B,l}^{(d)}(da)=\int_{\epsilon_0}^\infty
\mathcal{P}_{|a|}^{(r,|a|)}(\tau_{\epsilon_0}^{(Y)}>t) A^{d-1}\bar\mu_{B,l}^{(d)}(A)dA=\\
&\int_{\epsilon_0}^\infty\frac1{\mathcal{E}_A^{(r,A)}(\mathcal{U}_{r,A;\epsilon_0}(Y(t))|\tau_{\epsilon_0}^{(Y)}>t)}e^{-\lambda_0(r,A;\epsilon_0)\thinspace t}A^{d-1}\bar\mu_{B,l}^{(d)}(A)dA\sim\\
&\int_{\epsilon_0}^\infty\frac1{\mathcal{E}_A^{(r,A)}(\mathcal{U}_{r,A;\epsilon_0}(Y(t))|\tau_{\epsilon_0}^{(Y)}>t)}e^{-\lambda_0(r,A;\epsilon_0)\thinspace t}C(A)e^{-BA^l}dA,
\end{aligned}
\end{equation}
where
\begin{equation}\label{Casymp}
\lim_{A\to\infty}\frac{\log C(A)}{A^l}=0.
\end{equation}
The first equality in  \eqref{intasympBess} follows from
\eqref{sameprob1} and \eqref{sameprob2}, the second one follows from Theorem 3-Bessel, and the third one follows from
\eqref{targetdist}.
Since $\lambda_0(r,A;\epsilon_0)$ is decreasing to 0 as $A\to\infty$, it is clear that the asymptotic behavior of the right hand side of \eqref{intasympBess} as $t\to\infty$ depends only on large $A$.
Thus, in light of
\eqref{atoinfmult},
\begin{equation}\label{anotherasymp}
\begin{aligned}
&\int_{\epsilon_0}^\infty\frac1{\mathcal{E}_A^{(r,A)}(\mathcal{U}_{r,A;\epsilon_0}(Y(t))|\tau_{\epsilon_0}^{(Y)}>t)}e^{-\lambda_0(r,A;\epsilon_0)\thinspace t}C(A)e^{-BA^l}dA\sim\\
&\int_{\epsilon_0}^\infty C(A)e^{-\lambda_0(r,A;\epsilon_0)\thinspace t-BA^l}dA,\ \text{as}\ t\to\infty.
\end{aligned}
\end{equation}
Using \eqref{evequationmult}
and the fact that $K_{\frac{d-2}2}(x)\sim\sqrt{\frac\pi{2x}}e^{-x}$ as $x\to\infty$ \cite{AS,W},
we have
\begin{equation}\label{asympforeigen}
\lambda_0(r,A;\epsilon_0)\sim \eta\thinspace A^{\frac{1-d}2}e^{-\sqrt{(r-\lambda_0(r,A;\epsilon_0))\frac2D} A},\ \text{as}\ A\to\infty,\ \text{for some}\ \eta>0.
\end{equation}
Using \eqref{asympforeigen} and the fact that the $t\to\infty$ asymptotic behavior depends only on large $A$, we have for any $\delta>0$,
\begin{equation}\label{yetanother}
\begin{aligned}
&\int_{\epsilon_0}^\infty C(A)\exp(-t(\eta+\delta) A^{\frac{1-d}2}e^{-\sqrt{(r-\lambda_0(r,A;\epsilon_0))\frac2D} A}-BA^l)dA\le\\
&\int_{\epsilon_0}^\infty C(A)e^{-\lambda_0(r,A;\epsilon_0)\thinspace t-BA^l}dA\le\\
&\int_{\epsilon_0}^\infty C(A)\exp(-t(\eta-\delta) A^{\frac{1-d}2}e^{-\sqrt{(r-\lambda_0(r,A;\epsilon_0))\frac2D} A}-BA^l)dA,\\
&\text{for sufficiently large}\  t.
\end{aligned}
\end{equation}
By \eqref{Casymp} and the fact that
$\lim_{A\to\infty}\lambda_0(r,A;\epsilon_0)=0$, it follows that
for any $\delta>0$, we have for sufficiently large $A$,
\begin{equation}\label{twosidedestimate}
\begin{aligned}
&\exp(-t(\eta\pm\delta) e^{-\sqrt{(\frac{2r}D-\delta}) A}-(B+\delta)A^l)\le\\
&C(A)\exp(-t(\eta\pm\delta) A^{\frac{1-d}2}e^{-\sqrt{(r-\lambda_0(r,A;\epsilon_0))\frac2D} A}-BA^l)\le\\
&\exp(-t(\eta\pm\delta) e^{-\sqrt{(\frac{2r}D+\delta}) A}-(B-\delta)A^l).
\end{aligned}
\end{equation}
Applying \eqref{needed} to the integrals
$\int_{\epsilon}^\infty\exp(-t(\eta\pm\delta) e^{-\sqrt{(\frac{2r}D-\delta}) A}-(B+\delta)A^l)dA$ and
$\int_{\epsilon}^\infty\exp(-t(\eta\pm\delta) e^{-\sqrt{(\frac{2r}D+\delta}) A}-(B-\delta)A^l)dA$,
the proof of \eqref{asympprobtarget} now follows from
  \eqref{intasympBess}, \eqref{anotherasymp}, \eqref{yetanother} and \eqref{twosidedestimate}.
\hfill $\square$

\section{Proof of Proposition \ref{2dnoreset}}\label{2dno}
Recall that $P_x^{2;(0)}$ denotes probabilities for the standard two-dimensional Brownian motion starting from $x\in \mathbb{R}^2$.
By symmetry,  one has
\begin{equation}\label{2dbm}
P_0^{2;(0)}(\tau_a>t)=P_a^{2;(0)}(\tau_0>t).
\end{equation}
 This latter probability satisfies
\begin{equation}\label{2dformulahitting}
P_a^{2;(0)}(\tau_0>t)\approx1\wedge\frac{2\log \frac {|a|}{\epsilon_0}}{\log t},\ \text{for}\ t\ge2\ \text{and}\ |a|>\epsilon_0,
\end{equation}
 where $f(a,t)\approx g(a,t)$ means that there are constants $c_1,c_2>0$, independent of $a$ and $t$, such that
$c_1f(a,t)\le g(a,t)\le c_2f(a,t)$.
This follows from  the formula eight lines up from the bottom on page 774 of \cite{BMR}. The form there is slightly different from \eqref{2dformulahitting}, but is equivalent.
From \eqref{2dbm} and \eqref{2dformulahitting}, it follows that
\begin{equation}\label{deltapower}
\begin{aligned}
&\lim_{t\to\infty}P_0^{2;(0)}(\tau_{a_t}>t)=0,\ \text{if}\ \lim_{t\to\infty}\frac{|a_t|}{t^\delta}=0,\ \text{for all}\ \delta>0;\\
&\liminf_{t\to\infty}P_0^{2;(0)}(\tau_{a_t}>t)>0,\ \text{if}\ \lim_{t\to\infty}\frac{|a_t|}{t^\delta}>0,\ \text{for some}\ \delta>0.
\end{aligned}
\end{equation}
By Brownian scaling,
\begin{equation}\label{brownscale}
\lim_{t\to\infty}P_0^{2;(0)}(\max_{0\le s\le t}|X(t)|\ge |a_t|)=0,\ \text{if}\ \lim_{t\to\infty}\frac{|a_t|}{t^\frac12}=\infty.
\end{equation}
Now \eqref{2dnoresetresult}
follows from \eqref{deltapower} and \eqref{brownscale}.
\hfill $\square$

\section{Proof of Proposition \ref{speed}}\label{slowspeed}

We first consider the one-dimensional case. Clearly it suffices to consider the case that $a_t>0$.
By Theorem \ref{Main} and Proposition \ref{alimitunift},
it follows that
\begin{equation}\label{forspeedequ}
 \lim_{t\to\infty}P_0^{1;(r,0)}(\tau_{a_t}>t)=\begin{cases}0,\ \text{\rm if}\ \lim_{t\to\infty}t\lambda_0(r,0;a_t)=\infty;\\ 1,\ \text{\rm if}\  \lim_{t\to\infty}t\lambda_0(r,0;a_t)=0.\end{cases}
\end{equation}
Using \eqref{evequation}, we can replace $\lambda_0(r,0;a_t)$ in \eqref{forspeedequ}
by $re^{-a_t\thinspace\sqrt{\frac{2(r-\lambda_0(r,0;a_t))}D}}$. Thus,
\begin{equation}\label{forspeedequagain}
 \lim_{t\to\infty}P_0^{1;(r,0)}(\tau_{a_t}>t)=\begin{cases}0,\ \text{\rm if}\ \lim_{t\to\infty}\big(a_t-\sqrt{\frac D{2(r-\lambda_0(r,0;a_t))}}\thinspace\log t\big)=-\infty;\\ 1,\ \text{\rm if}\  \lim_{t\to\infty}\big(a_t- \sqrt{\frac D{2(r-\lambda_0(r,0,a_t))}}\thinspace\log t\big)=\infty.\end{cases}
\end{equation}
The case $\lim_{t\to\infty}(a_t-\sqrt{\frac D{2r}}\log t)=-\infty$ in
 \eqref{speedequ} follows from \eqref{forspeedequagain}.

 We now consider \eqref{speedequ} in the case that
  \begin{equation}\label{atgrowlog}
  \lim_{t\to\infty}(a_t-\sqrt{\frac D{2r}}\log t)=\infty.
  \end{equation}
    By Corollary \ref{cor1},
 $\lambda_0(r,0;a)\le re^{-c(r,a,D)\sqrt{\frac{2r}D}a}$, where $\lim_{a\to\infty}c(r,a,D)=1$.
Thus,
$$
\sqrt{\frac D{2(r-\lambda_0(r,0;a_t))}}\le
\sqrt{\frac D{2r}}\thinspace \big(1-e^{-c(r,a_t,D)\sqrt{\frac{2r}D}a_t}\big)^{-\frac12}.
%\le
%\sqrt{\frac D{2r}}\thinspace \big(1-e^{-c(r,a_t)\sqrt{\frac {2r}D}a_t}\big)^{-1}.
$$
Since \eqref{atgrowlog} holds, we have $e^{-c(r,a_t,D)\sqrt{\frac {2r}D}a_t}\le t^{-\frac12}$,
for all large $t$.
Consequently,
\begin{equation}\label{powerdecayt}
\sqrt{\frac D{2(r-\lambda_0(r,0;a_t))}}\le
\sqrt{\frac D{2r}}\thinspace (1+t^{-\frac12}),\ \text{for large}\ t.
\end{equation}
From \eqref{powerdecayt}, we conclude that if \eqref{atgrowlog} holds, then
$$
\lim_{t\to\infty}\big(a_t- \sqrt{\frac D{2(r-\lambda_0(r,0;a_t))}}\thinspace\log t\big)=\infty,
$$
and consequently, from \eqref{forspeedequagain},
$\lim_{t\to\infty}P_0^{1;(r,0)}(\tau_{a_t}>t)=1$. This concludes the proof of
\eqref{speedequ} in the case that \eqref{atgrowlog} holds.

We now turn to the multi-dimensional case. By Theorem 3-Bessel,  Proposition 4-Bessel,
\eqref{sameprob1} and \eqref{sameprob2},
it follows that
\begin{equation}\label{forspeedequmult}
 \lim_{t\to\infty}P_0^{d;(r,0)}(\tau_{a_t}>t)=\begin{cases}0,\ \text{\rm if}\ \lim_{t\to\infty}t\lambda_0(r,|a_t|;\epsilon_0)=\infty;\\ 1,\ \text{\rm if}\  \lim_{t\to\infty}t\lambda_0(r,|a_t|;\epsilon_0)=0.\end{cases}
\end{equation}
By \eqref{asympforeigen}, we can replace $\lambda_0(r,|a_t|;\epsilon_0)$ by
$|a_t|^{\frac{1-d}2}e^{-\sqrt{(r-\lambda_0(r,|a_t|;\epsilon_0))\frac2D} |a_t|}$ in \eqref{forspeedequmult}. Thus,
\begin{equation}\label{alternate1}
\begin{aligned}
&\lim_{t\to\infty}P_0^{d;(r,0)}(\tau_{|a_t|}>t)=0,\ \text{\rm if}\\
&\lim_{t\to\infty}\big(|a_t|-\sqrt{\frac D{2(r-\lambda_0(r,|a_t|;\epsilon_0))}}\log t+\frac{d-1}2
\sqrt{\frac D{2(r-\lambda_0(r,|a_t|;\epsilon_0))}}\log |a_t|\big)=-\infty;
\end{aligned}
\end{equation}
\begin{equation}\label{alternate2}
\begin{aligned}
& \lim_{t\to\infty}P_0^{d;(r,0)}(\tau_{a_t}>t)=1,\ \text{\rm if}\\
&\lim_{t\to\infty}\big(|a_t|-\sqrt{\frac D{2(r-\lambda_0(r,|a_t|;\epsilon_0))}}\log t+\frac{d-1}2
\sqrt{\frac D{2(r-\lambda_0(r,|a_t|;\epsilon_0))}}\log |a_t|\big)=\infty.
\end{aligned}
\end{equation}

We first consider the first case in
\eqref{speedequd} and thus assume that
\begin{equation}\label{conditioninthm1}
\lim_{t\to\infty}(|a_t|-\sqrt{\frac D{2r}}\log t+\gamma\log\log t)=-\infty,
\ \text{for
some}\ \gamma>\frac{d-1}2\sqrt{\frac D{2r}}.
\end{equation}
Without loss of generality we may assume that $\lim_{t\to\infty}|a_t|=\infty$, since otherwise
$\lim_{t\to\infty}P_0^{d;(r,0)}(\tau_{a_t}>t)=0$ follows trivially.
Then we have $\log |a_t|\le \log\log t+C$, for some constant $C$, and
 $\lim_{t\to\infty}\lambda_0(r,|a_t|;\epsilon_0)=0$.
From this and \eqref{conditioninthm1}, it follows that
the condition in the second line of \eqref{alternate1} holds, and consequently,
$\lim_{t\to\infty}P_0^{d;(r,0)}(\tau_{a_t}>t)=0$.

Now we consider the second case in \eqref{speedequd}
and thus assume that
\begin{equation}\label{conditioninthm2}
\lim_{t\to\infty}(|a_t|-\sqrt{\frac D{2r}}\log t+\frac{d-1}2\sqrt{\frac D{2r}}\log\log t)=\infty.
\end{equation}
By \eqref{asympforeigen}, we have for sufficiently large $A$,
$\lambda_0(r,A;\epsilon_0)\le re^{-\sqrt{(r-\lambda_0(r,A;\epsilon_0))\frac2D} A}$,
where we have included the $r$ for convenience in the next step.
Thus, for sufficiently large $t$,
\begin{equation}\label{negligible}
\begin{aligned}
&\sqrt{\frac D{2(r-\lambda_0(r,|a_t|;\epsilon_0))}}\le\sqrt{\frac D{2r}}(1-e^{-\sqrt\frac rD|a_t|})^{-\frac12}\le\\
&\sqrt{\frac D{2r}}(1-e^{-\sqrt\frac rD\sqrt\frac{D}{4r}\log t})^{-\frac12}=
\sqrt{\frac D{2r}}
(1-t^{-\frac12})^{-\frac12}\le\sqrt{\frac D{2r}}(1+t^{-\frac12}).
\end{aligned}
\end{equation}
From \eqref{negligible} and \eqref{conditioninthm2}, it follows that the  second line in \eqref{alternate2} holds,
and consequently,
$\lim_{t\to\infty}P_0^{d;(r,0)}(\tau_{B(A_t)}>t)=1$.
\hfill $\square$

\end{document}